 %plaintex

\magnification=1200

%caracteres gothiques
\font\tengoth=eufm10
\font\sevengoth=eufm7
\newfam\gothfam
\textfont\gothfam=\tengoth
\scriptfont\gothfam=\sevengoth
\def\goth{\fam\gothfam\tengoth}

%caracteres baton
\font\tenbboard=msbm10
\font\sevenbboard=msbm7
\newfam\bboardfam
\textfont\bboardfam=\tenbboard
\scriptfont\bboardfam=\sevenbboard
\def\bboard{\fam\bboardfam\tenbboard}

\newif\ifpagetitre
\newtoks\auteurcourant \auteurcourant={\hfill}
\newtoks\titrecourant \titrecourant={\hfill}

\pretolerance=500 \tolerance=1000 \brokenpenalty=5000
\newdimen\hmargehaute \hmargehaute=0cm
\newdimen\lpage \lpage=14.3cm
\newdimen\hpage \hpage=20cm
\newdimen\lmargeext \lmargeext=1cm
\hsize=11.25cm
\vsize=18cm
\parskip=0cm
\parindent=12pt

\def\margehaute{\vbox to \hmargehaute{\vss}}
\def\margebasse{\vss}

\output{\shipout\vbox to \hpage{\margehaute\nointerlineskip
  \corpsdepage\margebasse}
  \advancepageno \global\pagetitrefalse
  \ifnum\outputpenalty>-2000 \else\dosupereject\fi}

\def\corpsdepage{\hbox to \lpage{\hss\pagetexte\hskip\lmargeext}}
\def\pagetexte{\vbox{\makeheadline\pagebody\makefootline}}
\headline={\ifpagetitre\titleheadline \else
  \ifodd\pageno\rightheadline \else \leftheadline\fi\fi}
\def\leftheadline{\hfil\the\auteurcourant\hfil}
\def\rightheadline{\hfil\the\titrecourant\hfil}
\def\titleheadline{\hfill}
\pagetitretrue

%\font\petcap=amcsc10
\font\petcap=cmcsc10
\def\pc#1#2|{{\tenrm#1\sevenrm#2}}
\def\pd#1 {{\pc#1|}\ }

\def\pointir{\discretionary{.}{}{.\kern.35em---\kern.7em}\nobreak\hskip 0em
 plus .3em minus .4em }

\def\titre#1|{\message{#1}
             \par\vskip 30pt plus 24pt minus 3pt\penalty -1000
             \vskip 0pt plus -24pt minus 3pt\penalty -1000
             \centerline{\bf #1}
             \vskip 5pt
             \penalty 10000 }
\def\section#1|
               {\par\vskip .3cm
               {\bf #1}\pointir}
\def\ssection#1|
             {\par\vskip .2cm
             {\it #1}\pointir}

\long\def\th#1|#2\finth{\par\vskip 5pt
              {\petcap #1\pointir}{\it #2}\par\vskip 3pt}
\long\def\tha#1|#2\fintha{\par\vskip 5pt
               {\petcap #1.}\par\nobreak{\it #2}\par\vskip 3pt}

\def\rem#1|
{\par\vskip 5pt
                {{\it #1}\pointir}}
\def\rema#1|
{\par\vskip 5pt
                {{\it #1.}\par\nobreak }}

\def\proof{\noindent {\it Proof. }}

\def\qed{\quad\raise -2pt \hbox{\vrule\vbox to 10pt
{\hrule width 4pt
\vfill\hrule}\vrule}}

\def\cqfd{\ifmmode
\unkern\quad\hfill\qed
\else
\unskip\quad\hfill\qed\bigskip
\fi}

\newcount\n
\def\exo{\advance\n by 1 \par \vskip .3cm {\bf \the \n}. }

\def\ad{\mathop{\rm ad}\nolimits}

\def\Ad{\mathop{\rm Ad}\nolimits}

%produit semi-direct
\font\msbmten=msbm10
\def\ltimes{\mathbin{\hbox{\msbmten\char'156}}}

%%%%%%%%%%%%%%%%%%%%%%%%%%%%%%%%%%%%%

\hfill

\vskip 20 pt

\centerline{ANALYSIS OF  MINIMAL REPRESENTATIONS OF ${\rm SL}(n,{\bboard R})$  }

\vskip 20 pt
 
\centerline{Dehbia Achab}

\vskip 4 pt

\centerline{\tt dehbia.achab@imj-prg.fr}

\vskip 10 pt

\vskip 15pt

\noindent
{\bf Abstract}  Some minimal representations of ${\rm SL}(n,{\bboard R})$ can be realized on a Hilbert space of holomorphic functions. This is the analogue of the  Brylinski-Kostant model. They can also be realized on a Hilbert space of homogeneous functions on ${\bboard R}^n$. This is the analogue of the Kobayashi-Orsted model. We will describe the two realizations and  a transformation which maps one model to the other. It can be seen as an analogue of the classical Bargmann transform.

\bigskip
\noindent
Introduction
\vskip 2 pt
\noindent
1. The analogue of Brylinski-Kostant model, general case
\vskip 2 pt
\noindent
2.  Some harmonic analysis related to a quadratic form
\vskip 2 pt
\noindent
3. The Lie algebra construction, and isomorphism with ${\goth sl}(p+2,{\bboard C})$.
\vskip 2 pt
\noindent
4. The analogue of  the Kobayashi-Orsted model  for the minimal representations of  the group ${\rm SL}(p+2,{\bboard R})$.
\vskip 2 pt
\noindent
5. The analogue of the  Brylinski-Kostant model for  the minimal representations of  the group ${\rm SL}(p+2,{\bboard R})$.
\vskip 2 pt
\noindent
6. The intertwining operator
\vskip 2 pt
\noindent
References
\vskip 4 pt

\noindent
\bf Introduction. \rm  The construction of an analogue of the Bargmann-transform between Fock-type space and Schr\" odinger-type space of minimal representations  is given in a uniform manner,   in the case of simple real Lie groups of Hermitian type in [HKMO12].   We consider here a Lie group of  non Hermitian type, ${\rm SL}(n,{\bboard R})$,  and two minimal representations whose Fock-type spaces  and  Schr\" odinger-type spaces are described explicitely. We first recall  the theory  in [A11], [A12] and  [AF12]   to describe the analogue of the  Brylinski-Kostant model  ([BK94], [B97], [B98]),  in terms of holomorphic functions on a covering of the $K$-minimal nilpotent orbit.  The analogue of the Kobayashi-Orsted model ([KO03a], [KO03b], [KO03c]), in terms of homogeneous functions on ${\bboard R}^n$, will be given in  this paper.

\vskip 4 pt

 In [A11]   a general construction  for a simple  complex Lie  algebra ${\goth g}$ has been given, starting from a pair $(V,Q)$  where $V$ is a semi-simple Jordan algebra of rank $\leq 4$ and $Q$ a polynomial on $V$, homogeneous of degree 4.  The Lie algebra ${\goth g}$ is  of non Hermitian type. 
\vskip 1 pt
\noindent
In  [AF12],  the manifold $\Xi$ is the orbit of $Q$ under the conformal group ${\rm Conf}(V,Q)$  acting on a space of polynomials ${\cal W}$ by an irreducible representation ${\kappa}$.    The Fock space ${\cal F}(\Xi)$ is a Hilbert space  of holomorphic functions on the complex manifold $\Xi$. A spherical minimal representation of ${\goth g}$ is realized (when it exists) in ${\cal F}(\Xi)$. 
\vskip 1 pt
\noindent
In  [A12],   using the  decomposition $V=\oplus_{i=1}^sV_i$ of $V$  into simple summands, which means that $Q=\prod_{i=1}^s\Delta_i(z_i)^{k_i}$ gets factored as a product  of powers of Jordan determinants $\Delta_i$ of $V_i$,  and  considering irreducible representations $\kappa_i^{(k_i)}$  of   the conformal groups ${\rm Conf}(V_i,\Delta_i)$,  which act on   spaces ${\cal W}_i^{(k_i)}$ of the polynomials generated by the $\Delta_i^{k_i}(z_i-a_i)$ for $a_i\in V_i$,  the manifold $\Xi$ is  the orbit  of $Q$, under  the tensor product representation $\kappa^{(k_1,\ldots,k_s)}:=\otimes_{i=1}^s\kappa_i^{(k_i)}$ of the group   $\prod_{i=1}^s{\rm Conf}(V_i,\Delta_i)$ in ${\goth p}=\otimes_{i=1}^s{\cal W}_i^{(k_i)}$.  Then one gets a quotient map  $\widetilde\Xi\rightarrow \Xi, (\xi_1,\ldots,\xi_s) \mapsto \xi_1^{k_1}\ldots\xi_s^{k_s}$  of $\Xi$, where $\widetilde\Xi=\prod_{i=1}^s\Xi_i$ is the product of  the orbits $\Xi_i$ of  the  $\Delta_i$ under ${\rm Conf}(V_i,\Delta_i)$, acting  on  spaces ${\cal W}_i$ of the polynomials generated by the $\Delta_i(z_i-a_i)$  by irreducible representations $\kappa_i$.
The   ${\cal F}_{q}(\widetilde\Xi)$   
are  Hilbert spaces  of holomorphic functions on the complex manifold $\widetilde{\Xi}$, and  ${\cal F}_{q}(\widetilde\Xi)={\cal F}(\Xi)$  iff ${q}=0$. 
\vskip 4 pt
\noindent
Two minimal   representations of ${\goth g}$ are realized in the Fock-type spaces  ${\cal F}_{q}(\widetilde\Xi)$, for  some suitable multi-indices $ q\in {\bboard N}^s$. This is the analogue of the  Brylinski-Kostant model. 
\vskip 4 pt
\noindent
In this paper we consider the special case where $V={\bboard C}^{p}$,  $Q$ is the square of  a  quadratic form and the construction leads to the Lie algebra ${\goth sl}(p+2,{\bboard R})$ with $p \geq 3$. The two  minimal representations  $\rho_0$  and $\rho_1$ are   respectively realized   in ${\cal F}_0(\widetilde\Xi)$ and ${\cal F}_1(\widetilde\Xi)$.  
\vskip 4 pt
\noindent
On another hand, we realize   two minimal representations $\omega_0$ and $\omega_1$ of ${\rm SL}(p+2,{\bboard R})$  on  spaces  ${\cal V}_0({\bboard R}^{p+2})$ and  ${\cal V}_1({\bboard R}^{p+2})_1$ of homogeneous functions on ${\bboard R}^{p+2}$. This is the analogue of one of the Kobayashi-Orsted models.

\vskip 4 pt
\noindent
 We   also  give    explicit integral operators ${\cal B}_0$   from the space ${\cal V}_0({\bboard R}^{p+2})$ onto the  space ${\cal F}_0(\widetilde\Xi)$  which intertwines  the representations $\rho_0$ and $d\omega_0$ of ${\rm SL}(p+2,{\bboard R})$,  and ${\cal B}_1$ from the space ${\cal V}_1({\bboard R}^{p+2})$ onto the  space ${\cal F}_1(\widetilde\Xi)$ which intertwines  the representations $\rho_1$ and $d\omega_1$ of  ${\rm SL}(p+2,{\bboard R})$. These unitary operators are bijective and can be seen as analogues of the classical Bargmann transform.

\hfill
\eject

\noindent
\bf 1. The analogue of the Brylinski-Kostant model. General case \rm 
\vskip 10 pt

Let $V$ be  a  semi-simple Jordan algebra and $Q$  a homogeneous polynomial  on $V$. Let $L={\rm Str}(V,Q)$ be the structure group 

\vskip 4 pt

\centerline{${\rm Str}(V,Q)=\{g\in {\rm GL}(V) \mid \exists \gamma(g), Q(g\cdot z)=\gamma(g)Q(z)\}.$}

\vskip 4 pt
\noindent
 The conformal group ${\rm Conf}(V,Q)$ is the group of rational transformations $g$ of $V$ generated by: the translations $z\mapsto z+a$ ($a\in V$), the dilations $z\mapsto \ell \cdot z$
($\ell \in L$), and the inversion $\sigma :z\mapsto -z^{-1}$.  
\vskip 2 pt
Let $\cal W$ be the space of polynomials on $V$ 
generated by the translated $Q(z-a)$ of $Q$, with $a\in V$.  Let $\kappa$ be the cocycle representation of  ${\rm Conf}(V,Q)$ or of a covering of order two of it on $\cal W$,  defined in [A11] as follows:
\noindent
\vskip 1 pt
{\it Case 1.} In case there exists a character $\chi$ of Str$(V,Q)$ such that
$\chi^2=\gamma$, then let $K={\rm Conf}(V,Q)$. Define the cocycle 

\vskip 2 pt

\centerline{$\mu(g,z)=\chi((Dg(z)^{-1})\quad (g\in K,\ z\in V),$}

\vskip 2 pt
\noindent
and the representation $\kappa$ of $K$ on $\cal W$,
\vskip 2 pt

\centerline{$(\kappa(g)p)(z)=\mu(g^{-1},z)p(g^{-1}\cdot z).$}
\vskip 2 pt
\noindent
The cocycle $\mu(g,z)$ is a polynomial in $z$ of degree 
$\leq \, {\rm deg}\, Q$ and
\vskip 4 pt

\centerline{$(\kappa(\tau_a)p)(z)=p(z-a) \quad (a\in V),$}

\vskip 1 pt

\centerline{$(\kappa(\ell )p)(z)=\chi(\ell )p(\ell^{-1}\cdot z) \quad (\ell\in {\rm Str}(V,Q)),$}

\vskip 1 pt

\centerline{$(\kappa(\sigma)p)(z)=Q(z)p(-z^{-1}).$}

\vskip 4 pt
\noindent
{\it Case 2.} Otherwise the group $K$ is  a covering of ${\rm Conf}(V,Q)$, defined as the set of pairs
$\tilde g=(g,\mu)$, with  $g\in  {\rm Conf}(V,Q)$, and $\mu(\tilde g,\cdot)=\mu$ is a 
function on $V$ such that $\mu(z)^2=\gamma(Dg(z))^{-1}$, and 
$\bigl(\kappa(\tilde g)p\bigr)(z)=\mu({\tilde g}^{-1},z)p({\tilde g}^{-1}\cdot z)$.

\vskip 2 pt
We recall the construction of the Lie algebra ${\goth g}$  in [A.11]. The Jordan algebra is assumed to be of rank $\leq 4$ and the polynomial $Q$ is assumed to be homogeneous of degree 4.  We use the same notation 
${\goth k}={\rm Lie}(K), {\goth p}={\cal W}$. There is  $H\in {\goth l}={\rm Lie}(L)$ which defines  gradings of  ${\goth k}$ and ${\goth p}$:
\vskip 2 pt

\centerline{${\goth k}={\goth k}_{-1}+{\goth k}_0+{\goth k}_1,$}

\vskip 2 pt

\centerline{$ {\goth p}={\goth p}_{-2}+{\goth p}_{-1}+{\goth p}_0+{\goth p}_1+{\goth p}_2,$}

\vskip 1 pt
\noindent
with 

\centerline{${\goth k}_j=\{X\in {\goth k}\mid \ad (H)X=jX\}, \quad {\goth p}_j=\{p\in {\goth p} \mid {\rm d}\kappa(H)p=jp\},$}
 
 \vskip 2 pt
 
\centerline{${\goth k}_{-1}\simeq V,\quad {\goth k}_0=$Lie$(L),\quad
\Ad (\sigma): {\goth k}_j\rightarrow {\goth k}_{-j},$}

\vskip 2 pt

\centerline{${\goth p}_{-2}={\bboard C}, \quad {\goth p}_2={\bboard C}\, Q, 
\quad {\goth p}_{-1} \simeq V,\quad {\goth p}_1\simeq V, \quad \kappa(\sigma) : {\goth p}_j \rightarrow
{\goth p}_{-j}.$}

\vskip 2 pt
\noindent
Define   ${\goth g}={\goth k}\oplus {\goth p}$, $E:=Q$, $F:=1$. Then $[H,E]=2E, \quad [H,F]=-2F$ and (see Theorem 3.1 in [A11] and [AF12]) :
\vskip 2 pt
\bf Theorem 1.1. \rm 
There is a unique Lie algebra structure 
 on ${\goth g}$ satisfying 
 $$\eqalign{
&{\rm (i)}  \quad [X,X']=[X,X']_{\goth k} \quad (X,X' \in {\goth k}), \cr
&{\rm (ii)} \quad [X,p]={\rm d}\kappa(X)p \quad (X\in {\goth k}, p\in {\goth p}), \cr
&{\rm (iii)} \quad [E,F]=H. \cr}$$
\noindent
We recall also   the  real form  ${\goth g}_{\bboard R}$ of ${\goth g}$ which we consider. We fix a  Euclidean real form $V_{\bboard R}$ of the complex Jordan algebra $V$, denote by $z\mapsto
\bar z$ the conjugation of $V$ with respect to $V_{\bboard R}$, and then 
consider  the involution $g\mapsto \bar g$ of ${\rm Conf}(V,Q)$ given  by: 
$\bar g\cdot z=\overline{g\cdot \bar z}$.  
For $(g,\mu)\in K$ define
$\overline{(g,\mu)}=(\bar g,\bar\mu),\  {\rm where}\ 
\bar\mu(z)=\overline{\mu(\bar z)}.$
The  involution
$\alpha$ defined by $\alpha (g)=\sigma\circ\bar g\circ\sigma^{-1}$ is a Cartan
involution of $K$ (see [P.02] Proposition 1.1),  and  $K_{\bboard R}=\{g\in K\mid \alpha(g)=g\}$
is a compact real form of $K$.
For ${\goth u}$  the  compact real form of $\goth g$ such that ${\goth k}\cap {\goth u} = {\goth k}_{\bboard R}$ and  $\goth p_{\bboard R}:=\goth p\cap(i{\goth u})$, one  has   ${\goth p}_{\bboard R}=\{p\in {\goth p} \mid \kappa(\sigma)\bar p =p \}$, where $\bar p(z)=\overline{p(\bar z)}$, and  $\goth g_{\bboard R}:=\goth k_{\bboard R}+\goth p_{\bboard R}$  is a real form of $\goth g$ and the decomposition $\goth k_{\bboard R}+\goth p_{\bboard R} $ is its Cartan decomposition.

\vskip 2 pt

Using the decomposition of 
the Jordan algebra $V$ is a direct sum of simple ideals: 
$V=\sum_{i=1}^sV_i,$
and
$Q(z)=\prod _{i=1}^s\Delta_i(z_i)^{k_i}$, 
where $\Delta_i$ is the
determinant polynomial of the simple Jordan algebra 
$V_i$ and the $k_i$ are positive integers such that the degree of $Q$ is equal
to $\sum _{i=1}^sk_ir_i$, where $r_i$ is the rank of $V_i$., one can consider  the structure groups
$L_i={\rm Str}(V_i,\Delta_i)$,  the conformal groups ${\rm Conf}(V_i,\Delta_i)$,  the corresponding groups  $K_i$ analogous of $K$,  the spaces  
${\cal W}_i$  of the polynomials on $V_i$ 
generated by the   $\Delta_i(z_i-a_i)$ with $a_i\in V_i$  and the cocycle representations  
$\kappa_i$ of $K_i$, analogous of $\kappa$,  on the spaces  ${\cal W}_i$  given by $\bigl(\kappa_i(g_i)p\bigr)(z_i)=\mu_i(g_i^{-1},z_i)p(g_i^{-1}\cdot z_i),$ given,  
If $\chi_i$ exists,  by
$(\kappa_i(g_i)p)(z_i)=\mu_i(g_i^{-1},z_i)p(g_i^{-1}\cdot z_i)$and, when $\chi_i$ doesn't exist,  by  $\bigl(\kappa_i(\tilde g_i)p\bigr)(z_i)=\mu_i({\tilde g_i}^{-1},z_i)p({\tilde g_i}^{-1}\cdot z_i)$.

\vskip 2 pt
\noindent
Moreover, if ${\cal W}_i^{(k_i)}$ is   the  vector  space generated by the polynomials  $\Delta_i^{k_i}(z_i-a_i)$ for $a_i\in V_i$. Then, the  group $K_i$ acts on ${\cal W}_i^{(k_i)}$ by the representation  $\kappa_i^{(k_i)}$ given by
$(\kappa_i^{(k_i)}(g_i)p)(z_i)=\mu_i(g_i,z_i)^{k_i}p(g_i^{-1}\cdot z_i).$

\vskip 2 pt
\noindent
The product    $\widetilde L:=\prod_{i=1}^sL_i$ acting on $V$ by $g\cdot z=(g_i\cdot z_i)$ for  $g=(g_i)$,  $z=(z_i)$,  is a subgroup of $L$,  ${\rm Lie}(\widetilde L)={\rm Lie}(L)={\goth l}=\sum\limits_{i=1}^s{\goth l}_i$, $\gamma(g)=\prod_{i=1}^s\gamma_i^{k_i}(g_i).$
\vskip 2 pt
\noindent
The product  $\widetilde K:=\prod_{i=1}^sK_i $  and  $K$ have the same Lie algebra  $\goth k=\sum_{i=1}^s{\goth k}_i$, where 
${\goth k}_i=$Lie$(K_i)\simeq V_i\oplus{\goth l}_i\oplus V_i$.
\vskip 2 pt
\noindent
Furthermore,  the vector space ${\goth p}$ is  the tensor product of the ${\cal W}_i^{(k_i)}$  and the tensor product  representation $\kappa^{(k_1,\ldots,k_s)}:=\otimes_{i=1}^s\kappa_i^{(k_i)}$ of $\widetilde K$  in ${\goth p}$ and the representation $\kappa$ of $K$ in ${\goth p}$ have the same differential. Then,  the structure of simple Lie algebra  on ${\goth g}={\goth k}+{\goth p}$ can be obtained by considering  $\widetilde K$  (instead of  $K$) and $\kappa^{(k_1,\ldots,k_s)}$ (instead of $\kappa$).

\vskip 4 pt
The orbit $\Xi$ of $Q$ under the group ${\widetilde K}=\prod_{i=1}^sK_i$  acting on ${\goth p}$ by  $\kappa^{(k_1,\ldots,k_s)}$,  is  conical. It  is the $\widetilde K$-minimal nilpotent orbit in ${\goth p}$.
\vskip 4 pt
Let $\Xi_i=\{\kappa_i(g_i)\Delta_i\mid g_i\in  K_i\}$  be the
$K_i$-orbit  of $\Delta_i$ in ${\goth p}_i$. A polynomial $ \xi_i\in {\cal W}_i$ can be written
$\xi_i(v)=w_i\Delta_i(v)+\ {\rm terms \  of \ degree}\ < r_i \quad (w_i\in {\bboard C}),$
and $w_i=w_i(\nu_i)$ is a linear form on ${\cal W}_i$ which is semi-invariant under the preimage in $K_i$ of the  maximal 
parabolic subgroup $P_{\rm max}^{(i)}=L_i\ltimes N_i$, where $N_i$ is the group of translations $z_i\in V_i\mapsto z_i+a_i$, for $a_i\in V_i$.
The set $\Xi_{i,0}=\{\xi_i \in \Xi_i \mid w_i(\xi_i)\ne 0\}$
 is open and dense in $\Xi_i$.
 A polynomial $\xi_i\in \Xi_{i,0}$  can be written $\xi_i(v_i)=w_i\Delta_i(v_i-z_i) \quad
(w_i \in {\bboard C}^*, z_i\in V_i).$
Hence  we get a coordinate system $(w_i,z_i)\in {\bboard C}^*\times V_i$ for $\Xi_{i,0}$.
\vskip 4 pt
\noindent
In this coordinate  system, the cocycle action 
of $ K_i$ is given by 
\vskip 2 pt

\centerline{$\kappa_i(g_i) : (w_i,z_i) \mapsto \bigl(\mu_i(g_i,z_i)w_i,g_i\cdot z_i\bigr).$}

\vskip 4 pt
\noindent
Consider  $\widetilde\Xi=\prod_{i=1}^s\Xi_i$ and  the equivalence relation 
\vskip 2 pt

\centerline{$(\xi_1,\ldots,\xi_s) \sim (\lambda_1\xi_1,\ldots,\lambda_s\xi_s)$,  for $ (\lambda_1,\ldots,\lambda_s)\in {\bboard C}^s, \lambda_1^{k_1}\ldots\lambda_s^{k_s}=1$.}

\vskip 4 pt
\noindent
The  map   $\widetilde\Xi/\sim  \quad  \rightarrow \Xi$, $[(\xi_1,\ldots,\xi_s)] \mapsto \xi_1^{k_1}\ldots\xi_s^{k_s}$ is a diffeomorphism.

\vskip 2 pt
The group $K_i$  acts on the space  ${{\cal O}}(\Xi_i)$ of holomorphic functions     by: 

\vskip 2 pt

\centerline{$\bigl(\pi_i(g_i) f_i\bigr)({\xi_i})=f_i\bigl(\kappa_i(g_i^{-1})\xi _i\bigr)$.}

\vskip 2 pt
\noindent
If  $\xi_i(v_i)=w_i\Delta_i(v_i-z_i)$,  and $f_i\in {{\cal O}}(\Xi_i)$, we denote by
$ f_i(\xi_i)=\phi_i(w_i,z_i)$  the restriction of $f_i$ to $\Xi_{i,0}$. Then  $\phi_i\in {\cal O}({\bboard C}^*\times V_i)$. 
In the coordinates $(w_i,z_i)$, the  representation $\pi_i$ is given by
\vskip 4 pt

\centerline{$(\pi_i(g_i)\phi_i)(w_i,z_i)=\phi_i(\mu_i(g_i^{-1},z_i)w_i,g_i^{-1}\cdot z_i)$}
\vskip 4 pt
The group $\widetilde K$  acts on the tensor product  space $\otimes_{i=1}^s{{\cal O}}(\Xi_i)$   by the tensor product representation $\pi=\otimes_{i=1}^s\pi_i$, which is given by, for $g=(g_i)$:
\vskip 4 pt

\centerline{$\bigl(\pi(g)(f_1\otimes\ldots\otimes f_s)\bigr)(\xi_1,\ldots,\xi_s) =f_1(\kappa_1(g_1^{-1})\xi _1)\ldots f_s(\kappa_s(g_s^{-1})\xi _s)$}

\vskip 4 pt
\noindent
and in  coordinates,  by
\vskip 2 pt

\centerline{$(\pi(g)(\phi_1\otimes\ldots\otimes\phi_s))((w_1,z_1),\ldots,(w_s,z_s))=\prod_{i=1}^s\phi_i(\mu_i(g_i^{-1},z_i)w_i,g_i^{-1}\cdot z_i).$}

\vskip 2 pt
For  $m_i \in{\bboard Z}$,  the space  ${{\cal O}}_{m_i}(\Xi_i)$  of holomorphic functions
 $f_i$ on $\Xi_i$, homogeneous of degree $m_i$,   is invariant under the representation $\pi_i$.
If $f_i\in  {{\cal O}}_{m_i}(\Xi_i)$, then its 
restriction $\phi_i$ to $(\Xi_i)_0$ can be written 
$\phi_i(w_i,z_i)=w_i^{m_i}\psi_i(z_i)$
where $\psi_i$ is a holomorphic function on $V_i$. We  write 
$\widetilde {{\cal O}}_{m_i}(V_i)$ for the space of the functions $\psi_i$ corresponding to the
functions $f_i\in {{\cal O}}_{m_i}(\Xi_i)$, and denote by $\tilde\pi_{i,m_i}$ the representation of
$K_i$  on $\widetilde {{\cal O}}_{m_i}(V_i)$  corresponding to the restriction $\pi_{i,m_i}$ of
$\pi$ to ${{\cal O}}_{m_i}(\Xi_i)$. It is given by
\vskip 6pt

\centerline{$(\tilde\pi_{i,m_i}(g_i)\psi_i)(z_i)=\mu_i(g_i^{-1},z_i)^{m_i}\psi_i(g_i^{-1}\cdot z_i).$}

\vskip 2pt
\noindent
In particular 
\vskip 2 pt

\centerline{$(\tilde\pi_{i,m_i}(\sigma_i)\psi_i)(z_i)=\Delta_i^{m_i}(z_i)\psi_i(-z_i^{-1})$.}

\vskip 4 pt
\noindent
The spaces ${{\cal O}}_{m_i}(\Xi_i)$ ($=\{0\}$ if $m_i<0$), are 
  finite dimensional, and the
representations $\pi_{i,m_i}$ are   irreducible.

\vskip 4 pt
\noindent
The tensor product 
 ${{\cal O}}_{(m_1,\ldots,m_s)}(\widetilde\Xi)=\otimes_{i=1}^s{\cal O}_{m_i}(\Xi_i)$  is invariant under the   representation $\pi=\otimes_{i=1}^s\pi_i$. 
For  $f=f_1\otimes\ldots\otimes f_s\in {{\cal O}}_{(m_1,\ldots,m_s)}(\widetilde\Xi)$,  let  $\phi_i(w_i,z_i)$ be  the restriction of $f_i$ to $\Xi_{i,0}$,  then   $\phi_i(w_i,z_i)=w_i^{m_i}\psi_i(z_i)$ with $\psi_i\in  \widetilde{{\cal O}}_{m_i}(V_i)$, and $f$ is in coordinates given by 

\vskip 6 pt

\centerline{$(\phi_1\otimes\ldots\otimes\phi_s)((w_1,z_1),\ldots,(w_s,z_s))=w_1^{m_1}\ldots w_s^{m_s}\psi_1(z_1)\ldots\psi_s(z_s).$}

\vskip 4 pt
\noindent
For   $\widetilde {{\cal O}}_{(m_1,\ldots,m_s)}(V)=\otimes_{i=1}^s\widetilde{{\cal O}}_{m_i}(V_i)$ and  $\tilde\pi _{(m_1,\ldots,m_s)}=\otimes_{i=1}^s\tilde\pi_{m_i}$, 
\vskip 4 pt
\noindent

\centerline{$\bigl(\tilde\pi_{(m_1,\ldots,m_s)}(g)(\psi_1\otimes\ldots\otimes\psi_s) \bigr)(z_1,\ldots,z_s)=\prod\limits_{i=1}^s{\mu_i(g_i^{-1},z_i)}^{m_i}\psi_i(g_i^{-1}\cdot z_i).$}

\vskip 6 pt
We consider a Euclidean real form of each $V_i$, and the  compact real form $(K_i)_{\bboard R}$ of $K_i$ analogous of $K_{\bboard R}$. Then $\widetilde K_{\bboard R}:=\prod_{i=1}^s(K_i)_{\bboard R}$ is a compact real form of $\widetilde K$. There is a $(K_i)_{\bboard R}$-invariant norm on $\widetilde{\cal O}_{m_i}(V_i)$, which is unique up to a positive factor. It is given by 
$$\Vert \psi_i\Vert_{i,m_i}^2={1\over a_{i,m_i}}
\int_{V_i}\vert \psi_i(z_i)\vert^2H_i(z_i)^{-m_i}m_{i,0}(dz_i),$$
\noindent

\centerline{$a_{i,m_i}=\int_{V_i}H_i(z_i)^{-m_i}m_{i,0}(dz_i), $}

\vskip 4 pt
\noindent
and there  is a   $\widetilde K_{\bboard R}$-invariant norm  on ${{\cal O}}_{(m_1,\ldots,m_s)}(\widetilde\Xi)$, which is unique up to a positive factor, it is  given by
$\Vert  f\Vert_{(m_1,\ldots,m_s)}^2=\prod_{i=1}^s\Vert \psi_i\Vert_{i,m_i}^2$, 
where  $H_i(z_i)=H_i(z_i,z_i)$, with $H_i(z_i,z_i')$   holomorphic in $z_i$,
anti-holomorphic in $z_i'$ such that for $x_i \in(V_i)_{\bboard R}$,  $H_i(x_i,x_i)=\Delta_i(e_i+x_i^2)$,  
$m_{i,0}(dz_i)={1\over C_{i,0}}H_i(z_i)^{-2{n_i\over r_i}}m(dz_i),$ with

\vskip 4 pt
$$C_{i,0}=\int _{V_i}H_i(z_i)^{-2{n_i\over r_i}}m(dz_i),$$
\noindent
where   $n_i={\rm dim}(V_i)$and  the Lebesgue measure $m(dz_i)$ is  chosen such that $C_{i,0}=1$
(see [A12] Proposition 2.3 ).

\vskip 4 pt
First, the multi-index $q=(q_1,\ldots,q_s) \in {\bboard N}^s$ is determined  in such a way that the space of finite sums 
${\cal O}_{q}(\widetilde\Xi)_{\rm fin}$, 
which is $\widetilde K$-invariant,  carries a representation $\rho_{q}$ of the Lie algebra ${\goth g}$ (see [A.12a], Theorem 3.4).
For the existence of this representation one has to suppose that there is $i$ such that $q_i < k_i$, and to add  a  condition (T):
\vskip 4 pt

\centerline{$\exists \eta_q \quad {\rm such\quad  that} \quad \eta_{q}={q_i\over k_i}+{n_i\over k_ir_i} \quad (\forall i)\quad (T)$}
\vskip 2 pt
\noindent
For $X\in {\goth k}$,  $\rho_{q}(X)=d\pi(X)$.  In particular
\vskip 4 pt

\centerline{$\rho_q(H)=d\pi(H)={d\over dt}\mid_{t=0}\pi({\rm exp}tH).$}

\vskip 4 pt

\noindent
For $X\in {\goth p}$, $\rho_q(X)$ is a differential operator of degree $\leq 4$. The representation $\rho_{q}$ is determined by the operators $\rho_{q}(E)$ which involves the differential operator $Q({\partial\over\partial z})$, and $\rho_{q}(F)=-\rho_{q}(E)^*$.  More precisely, 
\vskip 4 pt

\centerline{$\rho_q(E)={\cal M}-\delta\circ{\cal D} \quad {\rm and} \quad \rho_q(F)={\cal M}^{\sigma}-\delta\circ{\cal D}^{\sigma}$}

\vskip 4 pt
\noindent
where ${\cal M}$ and ${\cal M}^{\sigma}$ are multiplication operators which map  the space ${\cal O}_m(\widetilde\Xi):={\cal O}_{(k_1m+q_1,\ldots,k_sm+q_s)}(\widetilde\Xi)$ into ${\cal O}_{m+1}(\widetilde\Xi)$, and  ${\cal D}$ and ${\cal D}^{\sigma}$ are differential operators which map  the space ${\cal O}_m(\widetilde\Xi)$ into ${\cal O}_{m-1}(\widetilde\Xi)$.  Their restrictions  to  ${\cal O}_m(\widetilde\Xi)$ are  :  for $\phi\in  {\cal O}_m(\widetilde\Xi)$, given by
\vskip 1 pt

\centerline{$\phi((w_1,z_1),\ldots,(w_s,z_s))=\prod_{i=1}^sw_i^{k_im+q_i}\psi_i(z_i),$}

\vskip 4 pt

\centerline{$({\cal M}\phi)((w_1,z_1),\ldots,(w_s,z_s))=\prod_{i=1}^sw_i^{k_i(m+1)+q_i}\psi_i(z_i)$}

\vskip 4 pt

\centerline{$({\cal D}\phi)((w_1,z_1),\ldots,(w_s,z_s))=\prod_{i=1}^sw_i^{k_i(m-1)+q_i}\bigl(\Delta_i^{k_i}\bigl({\partial\over\partial z_i}\bigr)\psi_i\bigr)(z_i),$}

\vskip 4 pt

\centerline{$({\cal M}^{\sigma}\phi)((w_1,z_1),\ldots,(w_s,z_s))=\prod_{i=1}^sw_i^{k_i(m+1)+q_i}\Delta_i^{k_i}(z_i)\psi_i(z_i)$}

\vskip 4 pt

\centerline{$({\cal D}^{\sigma}\phi)((w_1,z_1),\ldots,(w_s,z_s))\prod_{i=1}^sw_i^{k_i(m-1)+q_i}(({\bf D} _{-k_im+q_i}^{(i)})^*\psi_i)(z_i),$}

\vskip 4 pt
\noindent

\vskip 4 pt
\noindent
where ${\bf D}_{\alpha}^{(i)}$ is the Maass differential operator 
\vskip 2 pt

\centerline{${\bf D}_{\alpha }^{(i)}=\Delta_i(z_i) ^{k_i+\alpha }\Delta _i^{k_i}\Bigl({\partial \over \partial z_i}\Bigr)
\Delta _i(z_i)^{-\alpha},$}

\vskip 2 pt
\noindent
and  $({\bf D}_{\alpha }^{(i)})*$ is the adjoint of $({\bf D}_{\alpha }^{(i)})$ given by
$({\bf D}_{\alpha }^{(i)})*=J\circ {\bf D}_{\alpha }^{(i)}\circ J$, with $Jf(z_i)=f\circ \sigma_i (z_i)=f(-z_i^{-1})$,
and 
$\delta\phi=\delta_{(k_1m+q_1,\ldots,k_sm+q_s)}\phi$, 
 with the  sequence of numbers $(\delta_m)=(\delta_{( k_1m+q_1 ,\ldots, k_sm+q_s )})$   given, for $m\ne 0$,  by  
\vskip 4 pt

\centerline{$\delta_{(k_1m+q_1,\ldots,k_sm+q_s)}={1\over A(m+\eta_q )(m+\eta_q +1)}$}

\vskip 6 pt 
\noindent
and for $m=0$, by a case by case calculation (See [A12], Theorem 3.4).
\vskip 6 pt
Furthermore, we consider for a sequence $(c_m):=(c_{(k_1m+q_1,\ldots,k_sm+q_s)})$ of positive numbers, an inner product on ${\cal O}_q(\widetilde\Xi)_{\rm fin}$ such that, if $f=\sum_mf_m$, with $f_m\in {\cal O}_m(\widetilde\Xi)$,
$\Vert f\Vert^2=\sum_{m=0}^{\infty}{1\over c_m}\Vert f_m\Vert_m^2.$
This inner product is invariant under $K_{\bboard R}$. The completion of ${\cal O}_q(\widetilde\Xi)_{\rm fin}$ for this inner product is a Hilbert subspace ${\cal F}_q(\widetilde\Xi)$ of ${\cal O}_q(\widetilde\Xi)$. 

\vskip 2 pt
\noindent
The sequence $(c_m)=(c_{(k_1m+q_1,\ldots,k_sm+q_s)})$ is determined in such a way that the representation $\rho_q$ is unitary (see [A.12], Theorem 4.1). It is given, if $1-\eta_q$ is a root of the Bernstein polynomial, by
\vskip 4 pt

\centerline{$c_m=c_{(k_1m+q_1,\ldots,k_sm+q_s)}={(\eta_q +1)_m\over (\eta_q+a_0)_m(\eta_q +b_0)_m}{1\over m!}$}

\vskip 4 pt
\noindent
 and if not,  by
 \vskip 4 pt
 
 \centerline{$c_m=c_{(k_1m+q_1,\ldots,k_sm+q_s)}={(\eta_q+1 )_m(1)_m\over(\eta_q+a_0')_m (\eta_q +b_0')_m(\eta_q +c_0')_m}{1\over m!},$}
 
 \vskip 4 pt
\noindent
where the constants $a_0, b_0$ are given in Table 3, and $a_0',b_0'$ are given in Table 4   of that paper.

\noindent
\bf 2. Some harmonic analysis related to a quadratic form\rm  
\bigskip
\noindent
We consider on $V={\bboard C}^{p}$ the quadratic form 
\vskip 6 pt

\centerline{$\Delta(x)=x_1^2+\ldots+x_{p}^2.$}

\vskip 6 pt
\noindent
Then the structure group is 

\vskip 6 pt

\centerline{${\rm Str}(V,\Delta)={\bboard C}^*\times O(p,{\bboard C})$}

\vskip 6 pt
\noindent
(quotiented by $\{(1,I_p),(-1,-I_p)\}$). The conformal group  ${\rm Conf}(V,\Delta)$, generated by the translations, the structure group, and the inversion 
\vskip 6 pt

\centerline{$\sigma(z)=-{z\over \Delta(z)},$}

\vskip 6 pt
\noindent
is a complex Lie group isomorphic to $O(p+2,{\bboard C})$.  
\vskip 2 pt
The vector space ${\cal W}$ generated by the translated $\Delta(x-z)$  ($z\in V$) 

\vskip 2 pt

\centerline{$\Delta(x-z)=\Delta(x)-2\sum_{j=1}^pa_jx_j+\Delta(z),$}

\vskip 2 pt
\noindent
has dimension $p+2$. 
\vskip 2 pt
A transformation $l\in {\rm Str}(V,\Delta)$ is of the form $l\cdot z=\lambda u\cdot z$ with $\lambda\in {\bboard C}^*$, $u\in {\rm SO}(p,{\bboard C})$ and $\Delta(l\cdot z)=\lambda^2\Delta(z)$. Hence the characters $\gamma$ and $\chi$ of ${\rm Str}(V,\Delta)$ are given by 
$\gamma(l)=\lambda^2, \quad \chi(l)=\lambda$ 
and the cocycle $\mu$ by, for $g\in {\rm Conf}(V,\Delta)$ and $z\in V$, $\mu(g,z)=\chi(Dg(z)^{-1}).$

\vskip 4 pt
\noindent
The representation $\kappa$ of $K={\rm Conf}(V,\Delta)$ on ${\cal W}$ is defined by
\vskip 6 pt

\centerline{$(\kappa(g)\xi)(z)=\mu(g^{-1},z)\xi(g^{-1}\cdot z).$}

\vskip 4 pt

\noindent
In particular
\vskip 4 pt

\centerline{$(\kappa(\sigma)\xi)(z)=\Delta(z)\xi(-{z\over \Delta(z)}),$}

\vskip 4 pt

\centerline{
$(\kappa(\sigma)1)(z)=\Delta(z), \quad (\kappa(\sigma)\Delta)(z)=1.$}

\vskip 4pt
We consider the following basis ${\cal B}$ of ${\cal W}$:
\vskip 6 pt

\centerline{$e_j(x)=x_j  (1\leq j\leq p), e_{p+1}(x)={1\over 2}(\Delta(x)-1), e_{p+2}(x)={i\over 2}(\Delta(x)+1).$}

\vskip 6 pt
\noindent
We write an element $\xi\in {\cal W}$ as 
\vskip 6 pt

\centerline{$\xi(v)=\sum_{j=1}^{p+2}\xi_je_j(v).$}

\vskip 6 pt
\noindent
In particular 
\vskip 6 pt

\centerline{$\Delta=e_{p+1}-ie_{p+2}, \quad  1=-e_{p+1}-ie_{p+2}.$}
\vskip 6 pt
\noindent
Associated to this basis we  consider the bilinear form on ${\cal W}$ given by
$$\langle\xi,\eta\rangle=\sum_{j=1}^{p+2}\xi_j\eta_j.$$

\vskip 6 pt
\hfill
\eject
\noindent
For $g\in {\rm Conf}(V,\Delta)$ the matrix of $\kappa(g)$ with respect to this basis belongs to $O(p+2,{\bboard C})$, and this defines an isomorphism from ${\rm Conf}(V,\Delta)$ onto $O(p+2,{\bboard C})$. 
\vskip  2 pt
\noindent
The orbit $\widetilde\Xi$ of $\Delta$ under ${\rm Conf}(V,\Delta)$ is isomorphic to the isotropic  cone
\vskip 6 pt

\centerline{$\widetilde\Xi=\{\xi \in {\bboard C}^{p+2} \mid \sum_{j=1}^{p+2}\xi_j^2=0, \xi\ne 0\}.$}

\vskip 6 pt
\noindent
It is of dimension $p+1$, isomorphic to the homogeneous space 
\vskip 6 pt

\centerline{$O(p+2,{\bboard C})/O(p,{\bboard C})\ltimes {\bboard C}^p.$}

\vskip 6 pt
An element $\xi$ in ${\cal W}$ is a polynomial of degree $\leq 2$ which can be written

\vskip 6 pt

\centerline{$\xi(v)=w\Delta(v)+$ lower order  terms.}

\vskip 6 pt
\noindent
Hence, $\xi \mapsto w(\xi)$ is a linear form on ${\cal W}$. The set 
\vskip 6 pt

\centerline{$\widetilde\Xi_{0}=\{\xi \in \widetilde\Xi \mid w(\xi)\ne 0\}$}

\vskip 6 pt
\noindent
is open and dense in $\widetilde\Xi$, and 
\vskip 6 pt

\centerline{$w(\xi)={1\over 2}(\xi_{p+1}+i\xi_{p+2}).$}

\vskip 6 pt

An element $\xi\in \widetilde\Xi_{0}$ can be written 
$\xi(v)=w\Delta(v-z),$ 
with $w\in {\bboard C}^*$, $z\in V$, and the map
${\bboard C}^*\times V \rightarrow \widetilde\Xi_{0}, (w,z) \mapsto \xi(v)=w\Delta(v-z)$ 
is a diffeomorphism.  Explicitely
\vskip 6 pt

\centerline{$\xi(v)=w\Delta(v-z)=w(\Delta(v)-2\sum_{j=1}^pz_jv_j+\Delta(z))=\sum_{j=1}^{p+2}\xi_je_j(v),$}

\vskip 6 pt
\noindent
with
\vskip 6 pt

\centerline{$\xi_j=-2wz_j \quad (1\leq j\leq p),$}

\vskip 6 pt

\centerline{$\xi_{p+1}=w{1\over 2}(1-\Delta(z)),$}

\vskip 6 pt

\centerline{$\xi_{p+2}=w{1\over 2i}(1+\Delta(z)).$}

\vskip 6 pt
\noindent
One checks that 

\vskip 6 pt

\centerline{$\sum_{j=1}^{p+2}\xi_j^2=0,$}

\vskip 6 pt
\noindent
and also
\vskip 6 pt

\centerline{$\sum_{j=1}^{p+2}\Vert\xi_j\Vert^2={1\over 2}\vert w\vert^2H(z),$}

\vskip 6 pt
\noindent
with
\vskip 6 pt

\centerline{$H(z):=1+2\sum_{j=1}^p\vert z_j\vert^2+\vert\Delta(z)\vert^2.$}

\vskip 6 pt
\hfill
\eject
Recall  that,  for $m\geq 0$, ${\cal O}_m(\widetilde\Xi)$ denotes the space of holomorphic functions on $\widetilde\Xi$ which are homogeneous of degree $m$.  It is isomorphic to the quotient  ${\cal P}_m({\bboard C}^{p+2})/\Delta{\cal P}_{m-2}({\bboard C}^{p+2})$, where ${\cal P}_m({\bboard C}^{p+2})$ denotes the space of polynomials in $p+2$ variables,  homogeneous of degree $m$. 
\vskip 1 pt
The restriction to $\widetilde\Xi_{0}$ of a polynomial $f \in {\cal O}_m(\widetilde\Xi)$ can be written
\vskip 6 pt

\centerline{$u_m(\xi)=\phi_m(w,z)=w^m\psi(z)=\Bigr({\xi_{p+1}+i\xi_{p+2}\over 2}\Bigl)^m\psi_m(-{\xi_1\over 2w},\ldots,-{\xi_p\over 2w}),$}

\vskip 6 pt
\noindent
where $\xi$ is identified to $\xi=(\xi_1,\ldots,\xi_p,\xi_{p+1},\xi_{p+2})\in {\bboard C}^{p+2}$, and  $\psi_m$ is a holomorphic  polynomial on $V$ of degree $\leq m$. 
\vskip 2 pt
Observe that  the restriction of such polynomial $u_m$ to ${\bboard R}^{p+2}$  is an homogeneous harmonic polynomial of degree $m$ and the  holomorphic polynomial $\psi_m$ on ${\bboard C}^p$ is then given by 
\vskip 4 pt

\centerline{$\psi_m(z_1,\ldots,z_p)={1\over w^m}u_m(-2wz_1,\ldots,-2wz_p,{w\over 2}(1-\Delta(z)),{w\over 2i}(1+\Delta(z)))$}

\vskip 2 pt

\centerline{$=u_m(-2z_1,\ldots,-2z_p,{1-\Delta(z)\over 2},{1+\Delta(z)\over 2})$.}

\vskip 4 pt
\noindent
We get an isomorphism of ${\cal O}_m(\widetilde\Xi)$ to the space $\widetilde{\cal O}_m(V)$ of holomorphic polynomials $\psi_m$ on $V={\bboard C}^p$  such that  there exists an homogeneous harmonic polynomial $u_m$ of degree $m$ on ${\bboard R}^{p+2}$ such that
\vskip 6 pt

\centerline{$\psi_m(z_1,\ldots,z_p)=u_m(-2z_1,\ldots,-2z_p,{1-\Delta(z)\over 2},{1+\Delta(z)\over 2})$.}

\vskip 6 pt
The $O(p+2)$-invariant   Hilbert norm on ${\cal O}_m(\widetilde\Xi)$ is given, for   $f_m(\xi)=\phi_m(w,z)=w^m\psi_m(z)$,  by
$$\Vert f_m\Vert_m^2={1\over a_m}\int_V\vert\psi_m(z)\vert^2H(z)^{-m-p}m(dz),$$

\vskip 6 pt
\noindent
with
$$a_m=\int_VH(z)^{-m-p}m(dz).$$
Let ${\cal Y}_m({\bboard R}^{p+2})$ be the space of spherical harmonics of degree $m$ on ${\bboard R}^{p+2}$: harmonic polynomials which are homogeneous of degree $m$ on ${\bboard R}^{p+2}$.  Observe that, for $\xi\in \widetilde\Xi$, the function $x\mapsto \langle\xi,x\rangle^m$ belongs to ${\cal Y}_m$. 
The operator 
\vskip 6 pt

\centerline{${\cal A}_m : {\cal Y}_m({\bboard R}^{p+2}) \rightarrow {\cal O}_m(\widetilde\Xi), $}

$$({\cal A}_mu)(\xi)=\int_S\langle\xi,x\rangle^mu(x)s(dx),$$
\noindent
 intertwines the representations of $O(p+2)$  on both spaces. ($S$ is the unit sphere in ${\bboard R}^{p+2}$, and $s(dx)$ is the uniform  measure on $S$ with total measure equal to one). In fact, it follows from next Proposition and its corollary (see [F.15], Section 2).
\vskip 2 pt
\noindent
\bf Proposition 2.1| \rm For $u\in {\cal Y}_m$,
$$\int_S\langle\xi,x\rangle^mu(x)s(dx)=m!{\Gamma({p+2\over 2})\over\Gamma({p+2\over 2}+m)}u(\xi).$$

\vskip 6 pt
\noindent
\bf Corollary 2.2| \rm
 For $u\in {\cal Y}_m$, \quad $\Vert{\cal A}_mu\Vert_m^2=m!{\Gamma({p+2\over 2})\over\Gamma({p+2\over 2}+m)}\Vert u\Vert_{L^2(S)}^2.$

\vskip 8 pt

We consider now the case of the square of a  quadratic form: $V={\bboard C}^p$, and $Q(z)=\Delta^2(z)$, where  
\vskip 4 pt

\centerline{$\Delta(z)=z_{1}^2+\ldots+z_{p}^2.$}

\vskip 4 pt
\noindent
Then  
the vector space ${\cal W}^{(2)}$  generated by the translated $\Delta^2(x-z)$  with ($z\in V$) of the polynomial $Q=\Delta^2$,
\vskip 6 pt

\centerline{$\Delta^2(x-z)=(\Delta(x)-2\sum_{j=1}^{p}x_jz_j+\Delta(z))^2$}

\vskip 4 pt

\centerline{$=\Delta^2(x)+\Delta^2(z)+4\langle x,z\rangle^2+2\Delta(x)\Delta(z)-4\langle x,z\rangle\Delta(z)-4\langle x,z\rangle\Delta(x),$}

\vskip 6 pt
\noindent
has dimension  ${(p+4)(p+1)\over 2}$. We consider the following basis ${\cal B}^{(2)}$ of ${\cal W}^{(2)}$:
\vskip 4 pt

\centerline{$ e_0={1\over 2}(\Delta^2(x)-1), \quad \tilde e_0={i\over 2}(\Delta^2(x)+1),$}

\vskip 4 pt

\centerline{$e_j(x)=x_j^2 \quad (1\leq j\leq p),  \quad e_{jk}(x)=x_jx_k \quad (1\leq j<k\leq p),$}

\vskip 4 pt

\centerline{$e_j(x)=x_{j-p}\quad  (p+1\leq j\leq 2p), \quad e_j(x)=x_{j-2p}^3\quad  (2p+1\leq j\leq 3p).$}

\bigskip
\noindent
We write $\xi\in {\cal W}^{(2)}$ as 
\vskip 4 pt

\centerline{$\xi(x)=\sum\limits_{j=1}^{3p}\xi_je_j(x)+\sum\limits_{1\leq j<k\leq p}\xi_{jk}e_{jk}(x)+\xi_0e_0(x)+\tilde\xi_0\tilde e_0(x).$}

\vskip 4 pt
\noindent
For $\xi(x)=wQ(x-z)$, one has
\vskip 4 pt

\centerline{$\xi_0=w(1-Q(z)), \quad  \tilde\xi_0=w(-i-iQ(z)),$}

\vskip 4 pt

\centerline{$\xi_j= 2w\Delta(z)+4wz_j^2  \quad (1\leq j\leq p), \quad  \xi_{jk}=8wz_jz_k \quad (1\leq j<k\leq p),$}

\vskip 4 pt

\centerline{$\xi_j=-4w\Delta(z)z_{j-p} \quad  (p+1\leq j\leq 2p),    \xi_j=-4wz_{j -2p}\quad (2p+1\leq j\leq 3p).$}

\vskip 4 pt
\noindent
In particular,
\vskip 2 pt

\centerline{$Q=e_0-i\tilde e_0, \quad 1=-e_0-i\tilde e_0.$}

\vskip 6 pt
\noindent
  Denote by   $\Xi$  the orbit of $Q=\Delta^2$  under the action of the representation $\kappa^{(2)}$. Then  ${\widetilde\Xi}/\{\pm 1\}\simeq \Xi$. In coordinates,  the open and dense subset $\widetilde\Xi_0$ (resp.  $\Xi_0$)  of $\widetilde\Xi$ (resp. $\Xi$)  is diffeomorphic to  ${\bboard C}^*\times{\bboard C}^p$ (resp. $ ({\bboard C}^*/\{\pm 1\}) \times {\bboard C}^p$).
An element $\xi\in \Xi_0$ can be written 
$$\xi(v)=w\Delta^2(v-z)$$
\noindent
with
$$w=w(\xi)={1\over 2}(\xi_0+i\tilde\xi_0).$$

\vskip 4 pt
\noindent
One has ${\cal O}_m(\Xi)\simeq {\cal O}_{2m}(\widetilde\Xi)$,  and, formally, ${\cal O}_{m+{1\over 2}}(\Xi)\simeq {\cal O}_{2m+1}(\widetilde\Xi)$,    then  ${\cal F}_0(\widetilde\Xi)={\cal F}(\Xi)$ consists in functions on the orbit $\Xi$, but the space  ${\cal F}_1(\widetilde\Xi)$ does not.
\vskip 8 pt

\noindent
The  integral operators given by 
$${\cal A}_{2m} : {\cal Y}_{2m}({\bboard R}^{p+2})\rightarrow {\cal O}_{2m}(\widetilde\Xi),$$ 
$$({\cal A}_{2m}u)(\xi)=\int_{S}\langle\xi,x\rangle^{2m}u(x)s(dx),$$
\noindent
and
$${\cal A}_{2m+1} : {\cal Y}_{2m+1}({\bboard R}^{p+2})\rightarrow {\cal O}_{2m+1}(\widetilde\Xi),$$ 
$$({\cal A}_{2m+1}u)(\xi)=\int_{S}\langle\xi,x\rangle^{2m+1}u(x)s(dx),$$
are   isomorphisms.
In particular, for the functions $U_{2m} \in {\cal Y}_{2m}({\bboard R}^{p+2})$ and $U_{2m+1} \in {\cal Y}_{2m+1}({\bboard R}^{p+2})$
$$U_{2m}(x)=(x_{p+1}+ix_{p+2})^{2m},$$
\noindent
and
$$U_{2m+1}(x)=(x_{p+1}+ix_{p+2})^{2m+1},$$
\vskip 2 pt
\noindent
we get ${\cal A}_{2m}U_{2m}=\gamma(2m)F_{2m}$, and  ${\cal A}_{2m+1}U_{2m+1}=\gamma(2m+1)F_{2m+1}$ where
$$F_{2m}(\xi)=(\xi_{p+1}+i\xi_{p+2})^{2m},$$
\noindent
and
$$F_{2m+1}(\xi)=(\xi_{p+1}+i\xi_{p+2})^{2m+1},$$
\noindent
with
$$\gamma(2m)=2^{-2m}(2m)!{\Gamma({p\over 2}+1)\over\Gamma({p\over 2}+2+2m)}$$
\noindent
and
$$\gamma(2m+1)=2^{-2m-1}(2m+1)!{\Gamma({p\over 2}+1)\over\Gamma({p\over 2}+2+2m+1)}.$$
\vskip 4 pt
\noindent

\hfill
\eject
\noindent
\bf 3. The Lie algebra ${\goth g}$ and its isomorphism with ${\goth sl}(p+2,{\bboard C})$\rm 
\vskip 10 pt
\noindent
We describe the Lie algebra construction in the case  for $Q$ to be the square of a quadratic form, $Q=\Delta^2$. Let $H\in {\goth k}$ be the generator of one parameter group of dilations $h_t$ of $V$: $h_t\cdot z=e^{-t}z, \quad H={d\over dt}\mid_{t=0}h_t.$
\vskip 4 pt
\noindent
We get 
the matrix
\vskip 4 pt

\centerline{$\pmatrix{I_p&0&0\cr
0&\cosh t&i\sinh t\cr
0&-i\sinh t&\cosh t},$}

\vskip 4 pt
and  
\vskip 4 pt

\centerline{$d\kappa(H)=\pmatrix{0_p&0&0\cr
0&0&i\cr
0&-i&0}.$}
\vskip 4 pt
\noindent
We consider the elements $E,F\in {\goth p}:={\cal W}$ given by $E:=Q, F:=1$.
\vskip 2 pt
\noindent
By Theorem 1.1. there exists on ${\goth g}$ a unique Lie algebra structure such that for  $X,X' \in {\goth k}), p\in {\goth p})$,
\vskip 4 pt

\centerline{$ [X,X']=[X,X']_{\goth k} \quad [X,p]={\rm d}\kappa(X)p, \quad [E,F]=H.$}

\vskip 4 pt
\noindent
\bf Theorem 3.1| \rm If $Q=\Delta^2$ is the square of a quadratic form,
\vskip 4 pt

\centerline{$\Delta(z)=(z_1)^2+\ldots +(z_{p})^2,$}

\vskip 4 pt
\noindent
then  $({\goth g},{\goth k})$ is isomorphic to
$({\goth sl}(p+2,{\bboard C}),{\goth o}(p+2,{\bboard C}))$
\noindent 
 and  $({\goth g}_{\bboard R},{\goth k}_{\bboard R})$ is isomorphic to 
 $({\goth sl}(p+2),{\bboard R}),{\goth o}(p+2)).$
\vskip 2 pt
\noindent
\proof From the isomorphism ${\rm Conf}(V,\Delta) \rightarrow O(p+2,{\bboard C})$, 
we get  an isomorphism from the Lie algebra ${\goth k}$ of ${\rm Conf}(V,Q)$ onto  the Lie algebra 
${\goth o}(p+2,{\bboard C})$. Using the basis ${\cal B}^{(2)}$ of ${\cal W}^{(2)}$, we define an isomorphism from ${\goth p}={\cal W}^{(2)}$  onto 
$\{X\in {\rm Sym}(p+2,{\bboard C}), {\rm tr}(X)=0\}$, given by
\vskip 4 pt

\centerline{$\xi \mapsto \mu(\xi)=\pmatrix{\alpha(\xi)&\beta(\xi)\cr
\beta(\xi)^t&\delta(\xi) }$}

\vskip 4 pt
\noindent
with
\vskip 4 pt

\centerline{$\alpha(\xi)={1\over 4}\pmatrix{\xi_{11}&\xi_{12}&\ldots&\xi_{1p} \cr
\xi_{12}&\xi_{22}&\ldots&\xi_{2p}\cr
\vdots&\vdots&\vdots\cr
\xi_{1p}&\xi_{2p}&\ldots&\xi_{pp}}\quad   (\xi_{jj}=2\xi_j,   1\leq j\leq p),$}

\vskip 4 pt

\centerline{
$\beta(\xi)=-{1\over 4}\pmatrix{\xi_{p+1}&\xi_{2p+1}\cr
\vdots&\vdots\cr
\xi_{2p}&\xi_{3p}},$}

\vskip 4 pt

\centerline{$\delta(\xi)=\pmatrix{\xi_0&\tilde\xi_0\cr
\tilde\xi_0&-\xi_0-{\rm tr}(\alpha(\xi))}.$}

\vskip 4 pt
\noindent
 In the isomorphism described above, the images of the elements of the ${\goth sl}(2)$-triple $\{H,E,F\}$ are the following matrices
\vskip 4 pt

\centerline{$H\mapsto \widetilde H=\pmatrix{0_{pp}&&\cr
&0&i\cr
&-i&0},$}

\vskip 4 pt

\centerline{$E\mapsto  \widetilde E=\pmatrix{0_{pp}&&\cr
&1&-i\cr
&-i&-1},$}

\vskip 4 pt

\centerline{$F\mapsto  \widetilde F=\pmatrix{0_{pp}&&\cr
&1&i\cr
&i&-1}.$}

\vskip 4 pt
\noindent
One checks that $\{\tilde H,\tilde E,\tilde F\}$ is an ${\goth sl}(2)$-triple:
\vskip 4 pt

\centerline{$[\tilde H,\tilde E]=2\tilde E, \quad  [\tilde H,\tilde F]=-2\tilde F, \quad [\tilde E,\tilde F]=\tilde H.$}

\vskip 4 pt
\noindent
By using Theorem 1.1, this proves that we have obtained an explicit Lie algebra isomorphism from ${\goth g}$ onto ${\goth sl}(p+2,{\bboard C})$.

\vskip 8 pt
\noindent
One can notice that 
\vskip 4 pt

\centerline{${\goth k}=\{X\in {\goth sl}(p+2,{\bboard C}), X^t=-X\}$}

\vskip 4 pt
\noindent
and
\vskip 4 pt

\centerline{${\goth p}=\{X\in {\goth sl}(p+2,{\bboard C}), X^t=X\}.$}
\vskip 6 pt
\noindent
In particular, for $\xi(x)=Q(x-z)$, one has
\vskip 6 pt

\centerline{$\alpha(\xi)=2zz^t+\Delta(z)I_p =: a(z)\quad   (\xi_{jj}=2\xi_j,   1\leq j\leq p),$}

\vskip 6 pt

\centerline{$\beta(\xi)=\pmatrix{\Delta(z)z_1&z_1\cr
\vdots&\vdots\cr
\Delta(z)z_p&z_p}=: b(z),$}

\vskip 6 pt
\noindent
and  $\delta(\xi)=$
\vskip 6 pt

\centerline{$\pmatrix{1-Q(z)&-i-iQ(z)\cr
-i-iQ(z)&-1+Q(z)+(p-2)\Delta(z)}=:\pmatrix{d_0(z)&\tilde d_0(z)\cr
\tilde d_0(z)&-d_0(z)+(p-2)\Delta(z)}.$}

\vskip 6 pt
\noindent
Then the orbit $\Xi$ which we recall has dimension $p+1$,  consists in the matrices
\vskip 4 pt

\centerline{$m(z)=\pmatrix{a(z)&b(z)\cr
b(z)^t&d(z)}$}

\vskip 4 pt
\noindent
where 
\vskip 4 pt

\centerline{$d(z)=\pmatrix{d_0&-2i+id_0\cr
-2i+id_0&-d_0+(p-2)\Delta(z)},  \quad z\in {\bboard C}^p, d_0\in {\bboard C}.$}

\bigskip
\noindent
\bf 4. The analogue of  the Kobayashi-Orsted model  for the minimal representations of  the group ${\rm SL}(p+2,{\bboard R})$ \rm 
\vskip 10 pt
\noindent
Let $\Gamma$ be the  open  cone in ${\bboard R}^{p+2}$:
\vskip 6 pt

\centerline{$\Gamma=\{x \in {\bboard R}^{p+2} \mid \vert x\vert\ne 0\},$}

\vskip 6 pt
\noindent
and  $S$ be the unit sphere 
\vskip 6 pt

\centerline{$S=\{x\in {\bboard R}^{p+2} : \vert x\vert=1\}.$}

\vskip 6 pt
\noindent
The group $G_{\bboard R}={\rm SL}(p+2,{\bboard R})$ acts on ${\bboard R}^{p+2}$ by the natural representation, denoted by $L_g: x \mapsto gx$  ($g\in G_{\bboard R},x \in  {\bboard R}^{p+2}$). This action stabilizes the cone $\Gamma$. The multiplicative group ${\bboard R}_+^*$ acts on $\Gamma$ as a dilation and the quotient space $M=\Gamma/{\bboard R}_+^*$ is identified with $S$. This defines  an action of $G_{\bboard R}$ on $S$,  which leads to a $G_{\bboard R}$-equivariant principal ${\bboard R}_+^*$-bundle:
\vskip 6 pt

\centerline{$\Phi : \Gamma \rightarrow S, x\mapsto {x\over\vert x\vert}.$}

\vskip 6 pt
\noindent
\noindent
For $\lambda \in {\bboard C}$, let 
${\cal E}_{\lambda}(\Gamma)$ be  the space of ${\cal C}^{\infty}$-functions on $\Gamma$ homogeneous  of degree $\lambda$:
\vskip 6 pt

\centerline{${\cal E}_{\lambda}(\Gamma)=\{u\in {\cal C}^{\infty}(\Gamma) \mid u(tx)=t^{\lambda}u(x), \quad x\in \Gamma, t>0\}.$}

\vskip 6 pt
\noindent
The group ${\rm SL}(p+2,{\bboard R})$ acts naturally on ${\cal E}_{\lambda}(\Gamma)$, and, under  the action of the subgroup $O(p+2)$, the space   ${\cal E}_{\lambda}(\Gamma)$ decomposes as:
\vskip 6 pt

\centerline{${\cal E}_{\lambda}(\Gamma)\mid_S \simeq \bigoplus_{m=0}^{\infty}{\cal Y}_{m}({\bboard R}^{p+2}).$}

\vskip 6 pt
\noindent
Furthermore, for $\epsilon=\pm 1$, we put
\vskip 6 pt

\centerline{${\cal E}_{\lambda,\epsilon}(\Gamma):=\{u\in {\cal E}_{\lambda}(\Gamma) : u(-x)=\epsilon\cdot u(x), \quad x\in \Gamma\}.$}

\vskip 6 pt
\noindent
Then we have a direct sum decomposition into two  $G_{\bboard R}$-invariant subspaces:
\vskip 6 pt

\centerline{${\cal E}_{\lambda}(\Gamma)={\cal E}_{\lambda,1}(\Gamma)+{\cal E}_{\lambda,-1}(\Gamma)$}

\vskip 6 pt
\noindent
and, under the action of the subgroup $O(p+2)$, each of them  decomposes 
\vskip 6 pt

\centerline{${\cal E}_{\lambda,+1}(\Gamma)\mid_S \simeq \bigoplus_{m=0}^{\infty}{\cal Y}_{2m}({\bboard R}^{p+2})$}

\vskip 6 pt
\noindent 
and
\vskip 6 pt

\centerline{${\cal E}_{\lambda,-1}(\Gamma)\mid_S \simeq \bigoplus_{m=0}^{\infty}{\cal Y}_{2m+1}({\bboard R}^{p+2}).$}

\vskip 6 pt
\noindent
For $\lambda=\lambda_0=-(p+2)$, we denote  by:
\vskip 6 pt

\centerline{${\cal V}_0(\Gamma)={\cal E}_{\lambda_0,+1}, \quad {\cal V}_1(\Gamma)={\cal E}_{\lambda_0,-1}$}

\vskip 6 pt
\noindent
and by $\omega_0$ and $\omega_1$ the restrictions  to ${\cal V}_0(\Gamma)$ and ${\cal V}_1(\Gamma)$  of the natural action   
\vskip 6 pt

\centerline{$(\omega(g)u)(x):=u(g^{-1}x).$}
\bigskip
\noindent
\bf Lemma 4.1| \rm   For $u\in {\cal E}_{\lambda_0}(\Gamma)$,
$$-(d\omega)(\widetilde E)u)(x)=(x_{p+1}-ix_{p+2})\Bigl({\partial u\over\partial y_{p+1}}-i{\partial u\over\partial y_{p+2}}\Bigr),$$
$$-(d\omega)(\widetilde F)u)(x)=(x_{p+1}+ix_{p+2})\Bigl({\partial u\over\partial y_{p+1}}+i{\partial u\over\partial y_{p+2}}\Bigr).$$
\vskip 4 pt
\noindent
\bf Proof. \rm  For $M\in {\goth sl}(p+2,{\bboard C})$ and, for $u\in  {\cal E}_{\lambda_0}(\Gamma)$,
$$(d\omega(M)u)(x)={d\over dt}_{\mid_{t=0}}u({\rm exp}(tM)\cdot(x)=-Du(x)(Mx),$$
\noindent
where $D$ denotes the differential of the function $u$. For the matrix $M=F_0$ introduced in section 3 we get the formula of the lemma. \qed

\vskip 4 pt
\noindent
\bf Proposition  4.2| \rm   Consider the functions $U_{2m} \in  {\cal V}_0(\Gamma)$, $U_{2m+1}\in {\cal V}_1(\Gamma)$  defined by 
$$U_{2m}(x)={1\over\vert x\vert^{2m+p+2}}(x_{p+1}+ix_{p+2})^{2m},$$
$$U_{2m+1}(x)={1\over\vert x\vert^{2m+1+p+2}}(x_{p+1}+ix_{p+2})^{2m+1}.$$
\noindent
Then
$$d\omega_0(\widetilde E)U_{2m}=(2m-p-2)U_{2(m-1)},$$
$$d\omega_1(\widetilde E)U_{2m+1}=(2m+1-p-2)U_{2(m-1)+1}$$
\noindent
and
$$d\omega_0(\widetilde F)U_{2m}=(2m+p+2)U_{2(m+1)},$$
$$d\omega_1(\widetilde F)U_{2m+1}=(2m+1+p+2)U_{2(m+1)+1}.$$
\vskip 4 pt
\noindent
\bf Proof. \rm One writes
$$U_{\lambda}(x)={1\over\vert x\vert^{\lambda}}(x_{p+1}+ix_{p+2})^{\lambda}.$$
\noindent
and  uses the formula 
$$\Bigl({\partial \over\partial x_{p+1}}+i{\partial \over\partial x_{p+2}}\Bigr)(x_{p+1}+ix_{p+2})^{\lambda}=0,$$
\noindent
and
$$\Bigl({\partial \over\partial x_{p+1}}+i{\partial \over\partial x_{p+2}}\Bigr)\vert x\vert^{\lambda}=\lambda\vert x\vert^{\lambda-2}(x_{p+1}+ix_{p+2}).$$
\noindent

\vskip 4 pt
\noindent
\bf Theorem 4.3| \rm  
\vskip 2 pt
\noindent
1) (Irreducibility) The representations   $(\omega_0,{\cal V}_0(\Gamma))$  and $(\omega_1,{\cal V}_1(\Gamma))$  of ${\rm SL}(p+2,{\bboard R})$ are irreducible.
\vskip 2 pt
\noindent
2) ($K$-type decomposition) The underlying $({\goth g},K)$-modules,  $(\omega_0)_K$ and $(\omega_1)_K$, for $K=O(p+2,{\bboard C})$,  have the following  $K$-type formulas 
\vskip 4 pt

\centerline{
$(\omega_0)_K=\bigoplus\limits_{m=0}^{\infty}{\cal Y}_{2m}({\bboard R}^{p+2})$, \quad 
$(\omega_1)_K=\bigoplus\limits_{m=0}^{\infty}{\cal Y}_{2m+1}({\bboard R}^{p+2}).$}

\vskip 2 pt
\noindent
3) (Unitarity) The representations $\omega_0$ and $\omega_1$  of ${\rm SL}(p+2,{\bboard R})$ on ${\cal V}_0(\Gamma)$ and  respectively ${\cal V}_1(\Gamma)$ 
 are unitary  for the Hilbert norms on ${\cal V}_0(\Gamma)$ and ${\cal V}_1(\Gamma)$ defined as follows:  if  
\vskip 4 pt

\centerline{$u^{(0)}(x)=\sum_{m=0}^{\infty}{1\over  \vert x\vert^{2m+p+2}}u_{2m}(x), \quad u^{(1)}(x)=\sum_{m=0}^{\infty}{1\over\vert x\vert^{2m+1+p+2}}u_{2m+1}(x),$}

\vskip 4 pt
\noindent
 with  $u_{2m}\in {\cal Y}_{2m}({\bboard R}^{p+2})$,  $u_{2m+1}\in {\cal Y}_{2m+1}({\bboard R}^{p+2})$, then
$$\Vert u^{(0)}\Vert_{{\cal V}_0}^2=\sum\limits_{m=0}^{\infty}{1\over 2^{2m}}\int_{S}\vert u_{2m}(x)\vert^2s(dx),$$
$$ \Vert u^{(1)}\Vert_{{\cal V}_1}^2=\sum\limits_{m=0}^{\infty}{1\over 2^{2m}}\int_{S}\vert u_{2m+1}(x)\vert^2s(dx)$$
\vskip 2 pt
\noindent

This Theorem is less or more known : the representations $(\omega_0,{\cal V}_0(\Gamma))$  and $(\omega_1,{\cal V}_1(\Gamma))$ correspond respectively to the degenerate principal series $\pi_{\mu,0}^{{\rm GL}(p+2,{\bboard R})}$ and 
$\pi_{\mu,1}^{{\rm GL}(p+2,{\bboard R})}$ , for $\mu=0$, as they are given in [KOP11]  Section 4, in the case of ${\rm GL}(2n,{\bboard R})$.  
We give a proof  in our context.
\vskip 2pt
\noindent
\bf Proof. \rm 
1)   In fact, let ${\cal U}_0$ be an invariant subspace  of ${\cal V}_0$ for $\omega_0$  and ${\cal U}_1$ an invariant subspace  of ${\cal V}_1$ for $\omega_1$ and assume that ${\cal U}_0\ne\{0\}$ and ${\cal U}_1\ne\{0\}$. Since ${\cal U}_0$ and ${\cal U}_1$-invariant, and since the subspaces ${\cal Y}_{2m}({\bboard R}^{p+2})$ and ${\cal Y}_{2m+1}({\bboard R}^{p+2})$  are irreducible for  the restriction representations $(\omega_0)_K$ and $(\omega_1)_K$ respectively, then  there exist ${\cal I}_0, \subset {\bboard N}$ and ${\cal I}_1\subset {\bboard N}$ such that 
${\cal U}_0=\sum_{m\in {\cal I}_0} {\cal Y}_{2m}$ and ${\cal U}_1=\sum_{m\in {\cal I}_1} {\cal Y}_{2m+1}$. Furthermore,  since $d\omega_0(\tilde F)(U_{2m})=(2m+p+2)U_{2(m+1)}$ and $d\omega_1(\tilde F)(U_{2m+1})=(2m+p+3)U_{2(m+1)+1}$, then $d\omega_0(\tilde F)({\cal Y}_{2m})={\cal Y}_{2(m+1)}$ and $d\omega_1(\tilde F)({\cal Y}_{2m+1})={\cal Y}_{2(m+1)+1}$,  which means that  if  $m$ belongs to ${\cal I}_0$  then so does $m+1$ and if  $m'$ belongs to ${\cal I}_1$  then so does $m'+1$.  Moreover,  since $d\omega_0(\widetilde E)U_{2m}=(2m-p-2)U_{2(m-1)}$ and $d\omega_1(\widetilde E)U_{2m+1}=(2m+1-p-2)U_{2(m-1)+1}$ then $d\omega_0(\tilde E)({\cal Y}_{2m})={\cal Y}_{2(m-1)}$ and $d\omega_1(\tilde E)({\cal Y}_{2m+1})={\cal Y}_{2(m-1)+1}$,   which means that  if  $m$ belongs to ${\cal I}_0$  then so does $m-1$ and if  $m'$ belongs to ${\cal I}_1$  then so does $m'-1$.  It follows that ${\cal I}_0={\cal I}_1={\bboard N}$, i.e. ${\cal U}_0={\cal V}_0$  and ${\cal U}_1={\cal V}_1$. 
\vskip 2 pt
\noindent
2)  This follows from the decomposition to $K_{\bboard R}$-irreducible spaces:
\vskip 6 pt

\centerline{$L^2(S^{p+1}) \simeq \sum_{k=0}^{\infty}{\cal Y}_k({\bboard R}^{p+2}).$}

\vskip 6 pt
\noindent
3) Since the inner products on ${\cal V}_0$ and ${\cal V}_1$ are $K_{\bboard R}$-invariant,  then they are ${\goth sl}(p+2,{\bboard R})$-invariant if and only if, for every $X\in {\goth p}$, $d\omega(X)^*=-d\omega(\bar X)$. 
This is equivalent to the single condition
\vskip 4 pt

\centerline{$d\omega(\tilde E)^*=-d\omega(\tilde F).$}

\vskip 4 pt
\noindent
In fact, assume these conditions are satisfied. Then, for $X=g\tilde Eg^{-1}$ ($g\in K_{\bboard R}$),
\vskip 4 pt

\centerline{$d\omega(X)=\omega(g)d\omega(\tilde E)\omega(g)^{-1}, $}

\vskip 4 pt

\centerline{$d\omega(X)^*=-\omega(g^{-1})^*d\omega(\tilde F)\omega(g)^*=-\omega(g)d\omega(\tilde F)\omega(g)^{-1}=-d\omega(\bar X).$} 
\vskip 4 pt
\noindent
It remains to check that for  every $u, u' \in {\cal V}$, one has
\vskip 4 pt

\centerline{$\langle d\omega(\tilde E)u,u'\rangle_{\cal V}=-\langle u,d\omega(\tilde F)u'\rangle_{\cal V}$}
\vskip 4 pt
\noindent
which is equivalent to checking  that for every $u_k\in {\cal Y}_k, v_{k'}\in {\cal Y}_{k'}$, one has
\vskip 4 pt

\centerline{$\langle d\omega(\tilde E)u_k,v_{k'}\rangle_{\cal V}=-\langle u_k,d\omega(\tilde F)v_{k'}\rangle_{\cal V}.$}

\vskip 4 pt

For every $u_k\in {\cal Y}_k, v_{k'}\in {\cal Y}_{k'}$,  using integration by parts, one gets

$$\eqalignno{\langle &d\omega(\tilde E)u_{k},v_{k'}\rangle_{\cal V}=c_{k,k'}\int_S(x_{p+1}-ix_{p+2})({\partial u_k\over\partial x_{p+1}}-i{\partial u_k\over\partial x_{p+2}})\overline{v_{k'}(x)}s(dx)\cr
&=c_{k,k'}\int_S({\partial u_k\over\partial x_{p+1}}-i{\partial u_k\over\partial x_{p+2}})\overline{(x_{p+1}+ix_{p+2})v_{k'}(x)}s(dx)\cr
&=-c_{k,k'}\int_Su_k(x)\overline{(x_{p+1}+ix_{p+2})\Bigr({\partial v_{k'}\over\partial x_{p+1}}+i{\partial v_{k'}\over\partial x_{p+2}}}\Bigl)s(dx)\cr
&=-\langle u_{k},d\omega(\tilde F)v_{k'}\rangle_{\cal V}}$$

\noindent
where $c_{k,k'}$ is a constant.\quad $\qed$

\bigskip
Observe that 
the scalar products on ${\cal V}_0(\Gamma)$ and ${\cal V}_1(\Gamma)$ are given by
$$\langle u^{(0)},v^{(0)}\rangle_{{\cal V}_0}=\int_SD_p^{(0)}u^{(0)}(x)\overline{D_p^{(0)}v^{(0)}(x)}s(dx)$$
\noindent
and
$$\langle u^{(1)},v^{(1)}\rangle_{{\cal V}_1}=\int_SD_p^{(1)}u^{(1)}(x)\overline{D_pv^{(1)}(x)}s(dx)$$
\noindent
where $D_p^{(0)}$  and $D_p^{(1)}$ are the diagonal operators  given by
$$(D_p^{(0)}u^{(0)})(x)=\sum_{m=0}^{\infty}{1\over 2^m \vert x\vert^{2m+p+2}}u_{2m}(x).$$
\noindent
and
$$(D_p^{(1)}u^{(1)})(x)=\sum_{m=0}^{\infty}{1\over 2^m \vert x\vert^{2m+1+p+2}}u_{2m+1}(x).$$

\hfill
\eject
\noindent
\bf 5. The analogue of the  Brylinski-Kostant model for  the minimal representations of  the group ${\rm SL}(p+2,{\bboard R})$. \rm 
\vskip 10 pt
\noindent
We  assume $p\geq 3$. Following the method in [A12],  recalled in Section 1, we will construct the representations $\rho_0$  and $\rho_1$ on the spaces  of finite sums 
\vskip 4 pt

\centerline{${\cal F}_0(\widetilde\Xi)_{\rm fin}=\sum\limits_{m=0}^{\infty}{\cal O}_{2m}(\widetilde\Xi)$ and ${\cal F}_1(\widetilde\Xi)_{\rm fin}=\sum\limits_{m=0}^{\infty}{\cal O}_{2m+1}(\widetilde\Xi),$}

\vskip 4 pt
\noindent
More precisely,   we consider the case  where $V={\bboard C}^{p}$ with $p\geq 3$, $Q$ is the  square of a  quadratic form $\Delta$ : 
\vskip 4 pt

\centerline{$\Delta(v)=v_1^2+\ldots+v_{p}^2.$}

\vskip 4 pt
\noindent
 In this case, ${\goth g}={\goth sl}(p+2;{\bboard C})$, $K=O(p+2,{\bboard C})$, $s=1$, $k=2$, $r=2$.  The representations exist iff $q=0$ or $q=1$  (see [A12], Theorem 3.4 and Table 2). 
 In fact, condition (T) is fulfilled with $\eta_{q}={q\over 2}+{p+1\over 4}$ and, since one of the $q_i$ must be $< k=2$, then $q=0$ or $q=1$.
 \vskip 6 pt
The  Hilbert  spaces ${\cal F}_0(\widetilde\Xi)$ and ${\cal F}_1(\widetilde\Xi)$ decompose 
 \vskip 6 pt
 
\centerline{${\cal F}_0(\widetilde\Xi)=\bigoplus_{m=0}^{\infty}{\cal O}_{2m}(\widetilde\Xi),$ \quad ${\cal F}_1(\widetilde\Xi)=\bigoplus_{m=0}^{\infty}{\cal O}_{2m+1}(\widetilde\Xi),$}

\vskip 6pt

\noindent
where  the functions  $f_{2m}$ in  ${\cal O}_{2m}(\widetilde\Xi)$ and $f_{2m+1}$ in  ${\cal O}_{2m+1}(\widetilde\Xi)$ are, in coordinates, given by: 
\vskip 4 pt

\centerline{$\phi_{2m} (w,z)=w^{2m}\psi_{2m}(z) \quad (\psi_{2m}\in \widetilde{\cal O}_{2m}({\bboard C}^{p}), w \in {\bboard C}^*)$}

\vskip 2 pt

\centerline{$\phi_{2m+1} (w,z)=w^{2m+1}\psi_{2m+1}(z) \quad (\psi_{2m+1}\in \widetilde{\cal O}_{2m+1}({\bboard C}^{p}), w \in {\bboard C}^*).$}

\vskip 4 pt
\noindent
where   $\widetilde{\cal O}_{2m}({\bboard C}^{p})$ and $ \widetilde{\cal O}_{2m+1}({\bboard C}^{p})$ are respectively  the spaces of holomorphic polynomials $\psi_{2m}$ and $\psi_{2m+1}$ on ${\bboard C}^p$, such that there exist harmonic  homogeneous polynomials $u_{2m}$ of degree $2m$ and $u_{2m+1}$ of degree $2m+1$ on ${\bboard R}^{p+2}$ such that
\vskip 4 pt

\centerline{$\psi_{2m}(z_1,\ldots,z_p)=u_{2m}(-2z_1,\ldots,-2z_p,{1-\Delta(z)\over 2},{1+\Delta(z)\over 2})$,}

\vskip 2 pt

\centerline{$\psi_{2m+1}(z_1,\ldots,z_p)=u_{2m+1}(-2z_1,\ldots,-2z_p,{1-\Delta(z)\over 2},{1+\Delta(z)\over 2})$.}

\vskip 6 pt
The Euler operator ${\cal E}$ is defined as
\vskip 2 pt

\centerline{$({\cal E}\phi)(w,z)={d\over dt}\mid_{t=0}\phi(w,e^tz).$}

\vskip 2 pt
\noindent
One gets
\noindent
\vskip 4 pt

\centerline{$\rho_0(H)\phi_{2m}={\cal E}\phi_{2m}-2m\phi_{2m}$,}

\vskip 4 pt
\noindent

\centerline{ $\rho_1(H)\phi_{2m+1}={\cal E}\phi_{2m+1}-(2m+1)\phi_{2m+1}$.}

\vskip 6 pt
As in Section 1  we define  the operators ${\cal M}$ and ${\cal D}$ by
\vskip 4 pt

\centerline{$({\cal M}\phi)(w,z)=w^2\phi (w,z),$ and $({\cal D}\phi)(w,z)={1\over w^2}\bigr(Q({\partial\over\partial z})\phi\bigl) ((w,z))
$.}

\vskip 4 pt
\noindent
Hence ${\cal M}$ maps ${\cal O}_{2m}(\widetilde\Xi)$ into ${\cal O}_{2(m+1)}(\widetilde\Xi)$ and ${\cal D}$ into ${\cal O}_{2(m-1)}(\widetilde\Xi)$. 
\vskip 4 pt
Introduce the diagonal operator $\delta$: if $f=\sum_n f_n$ with $f_n\in {\cal O}_n(\widetilde\Xi)$, then
$\delta f=\sum\limits_n\delta_nf_n,$ 
where $(\delta_n)$ is a sequence of real numbers. 

\vskip 4 pt
\noindent
The operators  $\rho_0(F)$ and $\rho_0(E)$
 are   given by 
\vskip 4 pt

\centerline{$\rho_0(F)\phi_{2m} (w,z)=w^{2(m+1)}\psi(z)- \delta_{2(m-1)}w^{2(m-1)}\Delta^2({\partial\over\partial z})\psi(z),$}

\vskip 4 pt

\centerline{$\rho_0(E)\phi_{2m}(w,z)=w^{2(m+1)}\Delta^2(z)\psi(z) - \delta_{2(m-1)}w^{2(m-1)}{\bf D}_{-2m}^*\psi(z)$}

\vskip 4 pt
\noindent
and the  operators  $\rho_1(F)$ and $\rho_1(E)$ 
 are   given by 
\vskip 4 pt

\centerline{$\rho_1(F)\phi_{2m+1} (w,z)=w^{2(m+1)+1}\psi(z)- \delta_{2(m-1)+1}w^{2(m-1)+1}\Delta^2({\partial\over\partial z})\psi(z),$}

\vskip 4 pt

\centerline{$\rho_1(E)\phi_{2m+1}(w,z)=w^{2(m+1)}\Delta^2(z)\psi(z) - \delta_{2(m-1)+1}w^{2(m-1)+1}{\bf D}_{-2m-1}^*\psi(z)$.}

\vskip 4 pt
\noindent 
where   ${\bf D}_{\alpha}$ is  the Maass operator 
\vskip 2 pt

\centerline{${\bf D}_{\alpha}=\Delta(z)^{1+\alpha}\Delta({\partial\over\partial z})\Delta(z)^{-\alpha},$}

\vskip 4 pt

\centerline{$D^*=J\circ D\circ J$ with  $Jf(z)=f(-z^{-1}).$}
\vskip 2 pt
\noindent
and  the  sequences $(\delta_{2m})$ and $(\delta_{2m+1})$,   are  here given for $m\ne 0$  by 
\vskip 4 pt

\centerline{$\delta_{2m}={1\over 16(m+\eta_0)(m+\eta_0+1)}$ and $\delta_{2m+1}={1\over 16(m+\eta_1)(m+\eta_1+1)},$}

\vskip 4 pt
\noindent

\centerline{$\delta_0=1, \quad  \delta_1={1\over 4}.$}

\vskip 4 pt
\noindent
\bf Theorem   5.1| \rm Assume $p\geq  3$. Then $\rho_0$  and $\rho_1$ are  representations of the ${\goth sl}_2$-triple $\{E,F,H\}$
and extend as  irreducible representations of ${\goth g}$ on ${\cal F}_0(\widetilde\Xi)_{\rm fin}$ and  ${\cal F}_1(\widetilde\Xi)_{\rm fin}$ .

\vskip 2 pt
\noindent
This is a special case of  Theorem 3.4 and  Theorem 3.8 in [A12].
\vskip 8 pt
We consider for a sequence $(c_n)$ of positive numbers, an inner product on ${\cal F}_0(\widetilde\Xi)_{\rm fin}$ and an  inner product on ${\cal F}_1(\widetilde\Xi)_{\rm fin}$ such that, if $f=\sum_mf_{2m}$, with $f_{2m}\in {\cal O}_{2m}(\widetilde\Xi)$, and if $h=\sum_mh_{2m+1}$, with $h_{2m+1}\in {\cal O}_{2m+1}(\widetilde\Xi)$, 
\vskip 4 pt

\centerline{$\Vert f\Vert_{{\cal F}_0}^2=\sum\limits_{m=0}^{\infty}{1\over c_{2m}}\Vert f_{2m}\Vert_{2m}^2$ and $\Vert h\Vert_{{\cal F}_1}^2=\sum\limits_{m=0}^{\infty}{1\over c_{2m+1}}\Vert h_{2m+1}\Vert_{2m+1}^2.$}

\vskip 4 pt
\noindent
These inner products are  invariant under $K_{\bboard R}=SO(p+2)$. The completion  of ${\cal F}_0(\widetilde\Xi)_{\rm fin}$  for the first  inner product  and of ${\cal F}_1(\widetilde\Xi)_{\rm fin}$ for the second one  are  Hilbert subspaces ${\cal F}_0(\widetilde\Xi)$ and  ${\cal F}_1(\widetilde\Xi)$ of ${\cal O}(\widetilde\Xi)$ whose reproducing kernels are  given by
\vskip 2 pt

\centerline{${\cal K}_0(\xi,\eta)=\sum\limits_{m=0}^{\infty}c_{2m}\langle\xi,\bar{\eta}\rangle^{2m}$ and ${\cal K}_1(\xi,\eta)=\sum\limits_{m=0}^{\infty}c_{2m+1}\langle\xi,\bar{\eta}\rangle^{2m+1}$}

\vskip 4 pt
\noindent
\bf Theorem   5.2| \rm  We fix 
\vskip 6 pt

\centerline{$c_{2m}={(\eta_0+1)_m\over(\eta_0+{1\over 2})_m({1\over 2})_m}{1\over m!}$ and $c_{2m+1}={(\eta_1+1)_m\over(\eta_1+{1\over 2})_m({3\over 2})_m}{1\over m!}. $}

\vskip 6 pt
\noindent
Then, restricted to the real Lie algebra ${\goth g}_{\bboard R}$, $\rho_0$ is a unitary representation on the space ${\cal F}_0(\widetilde\Xi)_{\rm fin}$ and   $\rho_1$ is a unitary representation on  ${\cal F}_1(\widetilde\Xi)_{\rm fin}$.
\vskip 2 pt
\noindent
This is a special case of Theorem 4.1   in [A12].

\noindent
\bf 6. The  intertwining operator\rm 
\vskip 4 pt
\noindent
 Define the functions $b^{(0)}$ and  $b^{(1)}$ in one   complex variable by 
\vskip 6 pt

\centerline{$b^{(0)}(t)=\sum\limits_{m=1}^{\infty}b_{2m}t^{2m}$}

\vskip 6 pt
\noindent
and
\vskip 6 pt

\centerline{$b^{(1)}(t)=\sum\limits_{m=0}^{\infty}b_{2m+1}t^{2m+1}$}

\vskip 6 pt
\noindent
and the operators ${\cal B}_0: {\cal V}_0(\Gamma)_{\rm fin}\rightarrow  {\cal F}_0(\widetilde\Xi)_{\rm fin} $ and  ${\cal B}_1: {\cal V}_1(\Gamma)_{\rm fin}\rightarrow  {\cal F}_1(\widetilde\Xi)_{\rm fin} $:  for $x\in S$, $\xi \in \widetilde\Xi$,
$$({\cal B}_0u)(\xi)=\int_{S}b^{(0)}(\langle \xi,x\rangle)u(x)s(dx)$$
\noindent
and
$$({\cal B}_1u)(\xi)=\int_{S}b^{(1)}(\langle \xi,x\rangle)u(x)s(dx).$$
\noindent
By Proposition 2.1,  if the constants  $b_{2m}\ne 0$ and $b_{2m+1}\ne 0$,  then the restriction $({\cal B}_0)_{2m}=b_{2m}{\cal A}_{2m}$ of ${\cal B}_0$ to   ${\cal Y}_{2m}({\bboard R}^{p+2})$ and   the restriction  $({\cal B}_1)_{2m+1}=b_{2m+1}{\cal A}_{2m+1}$ of ${\cal B}_1$ to   ${\cal Y}_{2m+1}({\bboard R}^{p+2})$, are isomorphisms 
\vskip 4 pt

\centerline{$({\cal B}_0)_{2m}: {\cal Y}_{2m}({\bboard R}^{p+2})\rightarrow {\cal O}_{2m}(\Xi),$}

\vskip 4 pt

\centerline{$({\cal B}_1)_{2m+1}: {\cal Y}_{2m+1}({\bboard R}^{p+2})\rightarrow {\cal O}_{2m+1}(\Xi).
$}

\vskip 4 pt
\noindent
which intertwine the action of $O(p+2)$ by $\omega_0$ and $\rho_0$ and  respectively  $\omega_1$ and $\rho_1$  :  for $X\in {\goth k}$ with image $\widetilde X \in o(p+2,{\bboard C})$,

\vskip 4 pt

\centerline{${\cal B}_0d\omega_0(\tilde X)=\rho_0(X){\cal B}_0$ and  ${\cal B}_1d\omega_1(\tilde X)=\rho_1(X){\cal B}_1$.}

\vskip 4 pt
\noindent
\bf Theorem   6.1| \rm For 
\vskip 6 pt

\centerline{$b_{2m}={({p\over 2}+1)_{2m}\over 2^m({p\over 4}+{1\over 2})_m(2m)!}$ and 
$b_{2m+1}={({p\over 2}+1)_{2m+1}\over 2^m({p\over 4}+1)_m(2m+1)!} $}
 
 \vskip 6 pt
 \noindent
the operator  ${\cal B}_0$  
 intertwines the representations  $d\omega_0$ and $\rho_0$,  and the operator  ${\cal B}_1$  intertwines the representations $d\omega_1$ and $\rho_1$,  and they are  unitary and bijective.
\vskip 4 pt
\noindent
\bf Proof. \rm a) Consider the functions 
\vskip 6 pt

\centerline{$F_{2m}(\xi)=\Bigr({\xi_0+i\tilde\xi_0\over 2}\Bigl)^{2m}$ and $F_{2m+1}(\xi)=\Bigr({\xi_0+i\tilde\xi_0\over 2}\Bigl)^{2m+1}$}

\vskip 6 pt
\noindent
which correspond, in coordinates, to the functions
\vskip 6 pt

\centerline{$\Phi_{2m}(w,z)=w^{2m}$  and
$ \Phi_{2m+1}(w,z)=w^{2m+1}$.}

\vskip 6 pt
\noindent

\noindent
Then 
\vskip 6 pt

\centerline{$\rho_0(F)F_{2m}=F_{2(m+1)}$ and $\rho_1(F)F_{2m+1}=F_{2(m+1)+1}.$}

\vskip 10 pt
\noindent
Furthermore, 
$$\eqalign{&({\cal B}_0U_{2m})(\xi)=b_{2m}\int_{S}\langle \xi,x\rangle^{2m}(x_{p+1}+ix_{p+2})^{2m}s(dx)\cr
&=b_{2m}(2m)!{\Gamma\Bigr({p\over 2}+1\Bigl)\over\Gamma\Bigr({p\over 2}+1+2m\Bigl)}f_{2m}(\xi)\cr
&=\beta_{2m}F_{2m}(\xi).}$$
\noindent
and
$$\eqalign{&({\cal B}_1U_{2m+1})(\xi)=b_{2m+1}\int_{S}\langle \xi,x\rangle^{2m+1}(x_{p+1}+ix_{p+2})^{2m+1}s(dx)\cr
&=b_{2m+1}(2m+1)!{\Gamma\Bigr({p\over 2}+1\Bigl)\over\Gamma\Bigr({p\over 2}+2+2m\Bigl)}f_{2m+1}(\xi)\cr
&=\beta_{2m+1}F_{2m+1}(\xi)}$$
\noindent
with
\vskip 4 pt

\centerline{$\beta_{2m}=b_{2m}(2m)!{\Gamma\Bigr({p\over 2}+1\Bigl)\over\Gamma\Bigr({p\over 2}+1+2m\Bigl)}$ and $\beta_{2m+1}=b_{2m+1}(2m+1)!{\Gamma\Bigr({p\over 2}+1\Bigl)\over\Gamma\Bigr({p\over 2}+2+2m\Bigl)}.$}

\vskip 4 pt
\noindent
The intertwining relations 
\vskip 4 pt

\centerline{${\cal B}_0d\omega_0(\widetilde F)U_{2m}=\rho_0(F){\cal B}_0U_{2m}$ and ${\cal B}_1d\omega_1(\widetilde F)U_{2m+1}=\rho_1(F){\cal B}_1U_{2m+1}$}

\vskip 4 pt
\noindent
give  respectively  the conditions 
\vskip 4 pt

\centerline{$(2m+p+2)\beta_{2(m+1)}=\beta_{2m}$ and $(2m+1+p+2)\beta_{2(m+1)+1}=\beta_{2m+1}$}

\vskip 4 pt
\noindent
and if we fix $\beta_0=1$ and $\beta_1=1$, then
\vskip 6 pt

\centerline{$\beta_{2m}={1\over 2^m({p\over 4}+{1\over 2})_m}, \quad \beta_{2m+1}={1\over 2^m({p\over 4}+1)_m}.$}

\vskip 2 pt
\noindent
It follows that
\vskip 2 pt

\centerline{$b_{2m}={\Gamma\Bigr({p\over 2}+1+2m\Bigl)\over \Gamma\Bigr({p\over 2}+1\Bigl)}{1\over2^m({p\over 4}+{1\over 2})_m(2m)!}={({p\over 2}+1)_{2m}\over 2^m({p\over 4}+{1\over 2})_m(2m)!}$}

\vskip 2 pt
\noindent
and
\vskip 2 pt

\centerline{$b_{2m+1}={\Gamma\Bigr({p\over 2}+2+2m\Bigl)\over \Gamma\Bigr({p\over 2}+1\Bigl)}{1\over 2^m({p\over 4}+1)_m(2m+1)!}={({p\over 2}+1)_{2m+1}\over 2^m({p\over 4}+1)_m(2m+1)!} $.}

\vskip 6 pt
\noindent
One can also check the intertwining relations 

\vskip 4 pt

\centerline{${\cal B}_0d\omega_0(\widetilde E)U_{2m}=\rho_0(E){\cal B}_0U_{2m}$ and ${\cal B}_1d\omega_1(\widetilde E)U_{2m+1}=\rho_1(E){\cal B}_1U_{2m+1}$.}

\vskip 6 pt
Now,  recall that  ${\goth g}={\goth k}\oplus {\goth p}$ and ${\goth p}$ is a simple ${\goth k}$-module. Then ${\goth p}={\cal U}({\goth k})F$.    It follows that, since ${\cal B}_0$ intertwins the operators $d\omega_0(\widetilde F)$ and $\rho_0(F)$ and also  the  operators  $d\omega_0(\widetilde X)$ and $\rho_0(X)$ for every  $X\in {\goth k}$,  then  ${\cal B}_0$ intertwins the operators $d\omega_0(\widetilde Y)$ and $\rho_0(Y)$ for every $Y\in {\goth g}$.
Similarly, since ${\cal B}_1$ intertwins    the operators $d\omega_1(\widetilde F)$ and $\rho_1(F)$ and also  the  operators  $d\omega_1(\widetilde X)$ and  $\rho_1(X)$ for every  $X\in {\goth k}$,  then  ${\cal B}_1$ intertwins the operators $d\omega_1(\widetilde Y)$ and $\rho_1(Y)$ for every $Y\in {\goth g}$.
\vskip 2 pt
\noindent
b) We will check that, for $m\geq 0$,
\vskip 4 pt

\centerline{$\Vert{\cal B}_0U_{2m}\Vert_{{\cal F}}=\Vert U_{2m}\Vert_{{\cal V}_0}$ and $\Vert{\cal B}_1U_{2m+1}\Vert_{{\cal F}}=\Vert U_{2m+1}\Vert_{{\cal V}_1}.$}

\vskip 4 pt
\noindent
In fact, the norm of $u_{\lambda}(x)=(x_{p+1}+ix_{p+2})^{\lambda}$ in $L^2(S)$ is  given (see [F.15], Section 4) by the integral 
$$I_{\lambda}:=\int_S(x_{p+1}^2+x_{p+2}^2)^{\lambda}s(dx)=\lambda!{\Gamma({p\over 2}+1)\over\Gamma({p\over 2}+1+\lambda)}.$$
\noindent
Therefore,  the norm of $U_{2m}$ in ${\cal V}_0(\Gamma)$ and the norm of  $U_{2m+1}$ in ${\cal V}_1(\Gamma)$  are given by 
$$\eqalign{\Vert U_{2m}\Vert_{{\cal V}_0}^2&={1\over 2^{2m}}(2m)!{\Gamma({p\over 2}+1)\over\Gamma({p\over 2}+1+2m)}\cr
&={(2m)!\over 2^{2m}({p\over 2}+1)_{2m}}={(2m)!\over 2^{4m}({p\over 4}+{1\over 2})_m({p\over 4}+1)_m}}$$
\noindent
and
$$\eqalign{\Vert U_{2m+1}\Vert_{{\cal V}_1}^2&={1\over 2^{2m}}(2m+1)!{\Gamma({p\over 2}+1)\over\Gamma({p\over 2}+2+2m)}\cr
&={(2m+1)!\over 2^{2m}({p\over 2}+1)_{2m+1}}={(2m+1)!\over 2^{4m}({p\over 4}+1)_m({p\over 4}+{3\over 2})_m}}.$$
\noindent

\bigskip

The norm of $F_{2m}$ in the Hilbert space ${\cal F}_0(\widetilde\Xi)$ and  the norm of $F_{2m+1}$ in the Hilbert  space ${\cal F}_1(\widetilde\Xi)$ are given   by 
\vskip 4 pt

\centerline{$\Vert F_{2m}\Vert_{{\cal F}_0}^2={1\over c_{2m}}$ and $\Vert F_{2m+1}\Vert_{{\cal F}_1}^2={1\over c_{2m+1}}$}

\vskip 4 pt
\noindent
with
\vskip 4 pt

\centerline{$c_{2m}={({p\over 4}+1)_m\over({p\over 4}+{1\over 2})_m({1\over 2})_m}{1\over m!}$ and $c_{2m+1}={({p\over 4}+{3\over 2})_m\over({p\over 4}+1)_m({3\over 2})_m}{1\over m!}. $}

\vskip 4 pt
\noindent
In fact, the sequence $(c_m)$ is given in Theorem 4.1, table 3 and table 4  in [A12] 
\vskip 4 pt

\centerline{$c_{2m}={(\eta_0+1)_m\over(\eta_0+{1\over 2})_m({1\over 2})_m}{1\over m!}$ and $c_{2m+1}={(\eta_1+1)_m\over(\eta_1+{1\over 2})_m({3\over 2})_m}{1\over m!}. $}

\vskip 4 pt
\noindent
In our case, 
\vskip 4 pt

\centerline{$\eta_0={p\over 4}, \quad a_0={1\over 2}, \quad b_0={1\over 2}-{p\over 4}$}

\vskip 4 pt
\noindent
 and 
 \vskip 4 pt
 
 \centerline{$\eta_1={p\over 4}+{1\over 2}, \quad a_0=a_0'={1\over 2}, \quad b_0=b_0'={1\over 2}-{p\over 4}$}

\vskip 4 pt

\noindent
then 
\vskip 4 pt

\centerline{$c_{2m}={({p\over 4}+1)_m\over({p\over 4}+{1\over 2})_m({1\over 2})_m}{1\over m!}$ and 
$c_{2m+1}={({p\over 4}+{3\over 2})_m\over({p\over 4}+1)_m({3\over 2})_m}{1\over m!}. $}

\vskip 4 pt
\noindent
(one has to check that $1-\eta_0$ is a root of    the Bernstein polynomial but  $1-\eta_1$ isn't).
Therefore
$$\eqalign{\Vert{\cal B}_0U_{2m}\Vert_{{\cal F}_0}^2&=(\beta_{2m})^2\Vert f_{2m}\Vert_{{\cal F}_0}^2={\beta_{2m}^2\over c_{2m}}\cr
&={\beta^2({p\over 4}+{1\over 2})_m({1\over 2})_mm!\over 2^{2m}({p\over 4}+{1\over 2})_m^2({p\over 4}+1)_m}\cr
&={\beta^2({p\over 4}+{1\over 2})_m(2m)!\over 2^{4m}({p\over 4}+{1\over 2})_m^2({p\over 4}+1)_m}\cr
&={\beta^2(2m)!\over 2^{4m}({p\over 4}+{1\over 2})_m({p\over 4}+1)_m}}$$
\noindent
and
$$\eqalign{\Vert{\cal B}_1U_{2m+1}\Vert_{{\cal F}_1}^2&=(\beta_{2m+1})^2\Vert f_{2m+1}\Vert_{{\cal F}_1}^2={\beta_{2m+1}^2\over c_{2m+1}}\cr
&={{\beta'}^2({p\over 4}+1)_m({3\over 2})_mm!\over 2^{2m}({p\over 4}+1)_m^2({p\over 4}+{3\over 2})_m}\cr
&={{\beta'}^2({p\over 4}+1)_m(2m)!\over 2^{4m}({p\over 4}+1)_m^2({p\over 4}+{3\over 2})_m}\cr
&={{\beta'}^2(2m)!\over 2^{4m}({p\over 4}+1)_m({p\over 4}+{3\over 2})_m}}$$
\vskip 2 pt
\noindent
For $\beta=\beta'=1$,    one gets the unitarity of the operators ${\cal B}_0$ and ${\cal B}_1$.
\vskip 4 pt
\noindent
c)  Since the representations $d\omega_0$ and $d\omega_1$ are irreducible then the intertwining operators ${\cal B}_0$ and ${\cal B}_1$ are injective and since the representations $\rho_0$ and $\rho_1$ are irreducible then  ${\cal B}_0$ and ${\cal B}_1$ are surjective. 
\vskip 4 pt
\noindent
\bf Theorem   6.2| \rm For 
\vskip 6 pt

\centerline{$b_{2m}={({p\over 2}+1)_{2m}\over 2^m({p\over 2}+1)_m(2m)!}$ and 
$b_{2m+1}={({p\over 2})_{2m+1}\over 2^m({p\over 2}+{3\over 2})_m(2m+1)!}$}

\vskip 2 pt
the operators  ${\cal B}_0$ and ${\cal B}_1$ are unitary isomorphisms of $({\goth g},{\goth k})$ modules. The  irreducible unitary representation $\omega_0$ is equivalent to $\rho_0$ and the  irreducible unitary representation $\omega_1$ is equivalent to $\rho_1$.
\vskip 4 pt
\noindent

\centerline{\bf References}
\vskip 6pt
\noindent
\rm
\vskip 4 pt
\noindent
[A11] D. Achab (2011), \it Construction process for simple Lie algebras, \rm
J.  of Algebra, \bf 325,\rm 186-204.
\vskip 4 pt
\noindent
[AF12] D. Achab and J. Faraut (2012), \it Analysis of the Brylinski-Kostant model for minimal representations, Canad.J. of Math.,\bf 64, \rm721-754.
\vskip 4 pt
\noindent
[A12]  D. Achab (2012), \it  Minimal representations of simple real Lie groups of non Hermitian type, \rm arxiv.
\vskip 4 pt
\noindent
[AF16]  D. Achab and J. Faraut (2016), \it Analysis of the minimal representation of $O(m,n)$, \rm Preprint.
\vskip 4 pt
\noindent
[B97] R. Brylinski (1997), \it Quantization of the 
4-dimensional nilpotent orbit of ${\rm SL}(3,{\bboard R})$, \rm 
Canad.  J. of Math., \bf 49, \rm 916-943.
\vskip 4 pt
\noindent
[B98]  R. Brylinski (1998), \it Geometric quantization of real minimal nilpotent orbits, \rm Symp. Geom.,
Diff. Geom.  and Appl., \rm 9, 5-58.
\vskip 4pt
\noindent
[BK94] R. Brylinski and B. Kostant (1994), \it 
Minimal representations, geometric quantization and unitarity, \rm 
Proc. Nat. Acad., \bf 91, \rm 6026-6029.
\vskip 4 pt
\noindent
[F15] J. Faraut (2015), \it Analysis of the minimal representation of $O(n,n)$, \rm Preprint.
\vskip 4 pt
\noindent
[HKMO12]  (2012) J. Hilgert, T. Kobayashi,  J. M\"ollers, B. Orsted (2012), \it  Fock model and Segal-Bargmann transform  for minimal representations of Hermitian Lie groups, \rm Journal of Functional Analysis,\bf  263, \rm 3492-3563.
\vskip 4 pt
\noindent
[KO03a] T. Kobayashi and B. Orsted (2003), \it Analysis on the minimal representation of $O(p,q)$, I. Realization via conformal geometry, \rm Adv. Math., \bf 180, \rm  486-512.
\vskip 4 pt
\noindent
[KO03b] T. Kobayashi and B. Orsted (2003), \it Analysis on the minimal representation of $O(p,q)$, II. Branching laws, \rm Adv. Math.,\bf 180, \rm 513-550.
[KO03c] T. Kobayashi and B. Orsted (2003), \it Analysis on the minimal representation of $O(p,q)$, III.  Ultrahyperbolic equations on ${\bboard R}^{p-1,q-1}$, \rm Adv. Math., \bf 180, \rm 551-595.
\vskip 4 pt
\noindent
[KM11] T. Kobayashi and G. Mano (2011), \it The Schr\"odinger model for the minimal representation of the indefinite orthogonal group $O(p,q)$, \rm Memoir of the American Mathematical Society, Volume 213, Number 1000.
\vskip 4 pt
\noindent
[KOP11]  T. Kobayashi, B. Orsted, and M. Pevzner (2011), \it Geometric analysis on small unitary representations of ${\rm GL}(n,{\bboard R})$, \rm J. Funct. Anal. \bf 260 \rm,  1682-1720.
\vskip 4 pt
\noindent
[M] V.F. Molchanov (1970), \it Representations of pseudo-orthogonal groups associated with a cone, \rm Math. USSR Sbornik, \bf 10, \rm 333-347.

\end{document}